\documentclass[onefignum,onetabnum]{siamart190516} % arXiv

\usepackage{amsfonts}
\usepackage[utf8]{inputenc}
\usepackage{graphicx}
\usepackage{epstopdf}
\usepackage{algorithmic}
\usepackage{xcolor}
\usepackage{pgfplots}
\pgfplotsset{compat=1.15}
\usepackage{xfrac}
\usepackage{mathtools}
\usepackage{booktabs}
\usepackage[style=base]{caption}
\usepackage{pgfplotstable}
\ifpdf
  \DeclareGraphicsExtensions{.eps,.pdf,.png,.jpg}
\else
  \DeclareGraphicsExtensions{.eps}
\fi

\usepackage{geometry,marginnote}

\allowdisplaybreaks[4]

\newsiamremark{remark}{Remark}
\newsiamremark{hypothesis}{Hypothesis}
\crefname{hypothesis}{Hypothesis}{Hypotheses}
\newsiamthm{claim}{Claim}

\headers{Oscillations in a Becker--D\"oring model}{B. Niethammer, R.L. Pego, A. Schlichting, and J. J. L. Vel\'{a}zquez}

\title{Oscillations in a Becker--D\"oring model\\ with injection and depletion
\thanks{Submitted to the editors \today.
\funding{The authors gratefully acknowledge the support of the Hausdorff Research Institute for Mathematics (Bonn), through the Junior Trimester Program on \emph{Kinetic Theory}. 
This work is supported  by the German Research Foundation (DFG) under Germany's Excellence Strategy EXC 2044 -- 390685587,
\emph{Mathematics M\"unster: Dynamics--Geometry--Structure} (AS), EXC-2047/1 -- 390685813, the \emph{Hausdorff Center for Mathematics}, as well as the Collaborative Research Center 1060 -- 211504053, \emph{The Mathematics of Emergent Effects} at the University of Bonn (BN, AS, JV). 
This material is also based upon work supported by the National Science Foundation under Grant Nos.\ DMS 1812609 and 2106534 (RLP).}}}

\author{B. Niethammer\thanks{Institute of Applied Mathematics, University Bonn
  (\email{niethammer@iam.uni-bonn.de}, \email{velazquez@iam.uni-bonn.de}).}
\and R.L. Pego\thanks{Department of Mathematical Sciences, Carnegie Mellon University
  (\email{rpego@cmu.edu}).}
\and A. Schlichting\thanks{Applied Mathematics M\"unster, %Institute for Analysis and Numerics,
	 University M\"unster (\email{a.schlichting@uni-muenster.de})}
\and J. J. L. Velázquez\footnotemark[2]}

\usepackage{amsopn}

\DeclareMathOperator{\re}{Re}
\DeclareMathOperator{\im}{Im}
\DeclareMathOperator{\sign}{sign}

\DeclareMathOperator{\nucl}{nuc}
\DeclareMathOperator{\cond}{con}
\DeclareMathOperator{\transport}{transport}

\newcommand{\C}{\mathbb C}
\newcommand{\N}{\mathbb N}
\newcommand{\R}{\mathbb R}
\newcommand{\Z}{\mathbb Z}

\newcommand{\eps}{\varepsilon}

\newcommand{\crit}{\text{crit}}
\newcommand{\rem}{\text{rm}}

\newcommand{\ignore}[1]{}
\newcommand{\HIDE}[1]{}  % {#1} to show
\newcommand{\vth}{\vartheta}  % rescaled theta

\newcommand{\rfac}{\eta}
\newcommand{\rexp}{r}

\renewtheorem*{lemma}[theorem]{Lemma}  % Removes the number. There is only one lemma, section 3.1
\renewtheorem*{proposition}[theorem]{Proposition}  

\newcommand{\AppendixSupp}{ the Appendix} % arXiv

\begin{document}

\maketitle

\begin{abstract}
  We study the Becker--Döring bubblelator, 
  a variant of the Becker--Döring coagulation-fragmenta\-tion system that models the growth of clusters by gain or loss of monomers. 
  Motivated by models of gas evolution oscillators from physical chemistry, 
  we incorporate injection of monomers and depletion of large clusters.
  For a wide range of physical rates, the Becker--Döring system itself exhibits a dynamic phase transition as mass density increases past a critical value.
  We connect the Becker--Döring bubblelator to a transport equation coupled with an integrodifferential equation for the excess monomer density by formal asymptotics in the near-critical regime.
  For suitable injection/depletion rates, we argue that time-periodic solutions appear via a Hopf bifurcation.  
  Numerics confirm that the generation and removal of large clusters can become desynchronized, 
  leading to temporal oscillations associated with bursts of large-cluster nucleation.
\end{abstract}



\section{Introduction}\label{S.intro}

Becker and Döring~\cite{BD1935} provided one of the original descriptions of a mechanism of particle growth in the theory of nucleation from supersaturated vapor. 
The main assumptions of their model are that individual clusters consist of atomic parts called \emph{monomers}, 
and that the growth and shrinkage of clusters occurs only by the addition and removal of single monomers. 
Although this process is not necessarily realized by chemical kinetics, it is convenient to be interpreted as a reaction network of the form
\begin{equation}\label{e:BD:ChemReact}
     \{1\} + \{k\} \xrightleftharpoons[b_{k+1}]{a_k} \{k+1\} \:,  \qquad k=1,2,3,\dots \, .
\end{equation}
As noted by Slemrod~\cite{slemrod2000}, the Becker-D\"oring equations provide perhaps the simplest model capable of a realistic description of several phenomena associated with the dynamics of phase changes.
Starting from the seminal work of Ball, Carr and Penrose~\cite{BCP86}, the mathematical theory for these equations has been developed in great detail.
Many aspects of the long-time behavior of solutions and the implications for the emergence of phase transitions are understood, 
but there are also still open questions; see~\cite{HingantYvinec2017} for a recent review.

In this work we add to~\eqref{e:BD:ChemReact} two reaction mechanisms, which are motivated by the dynamics of 
chemical oscillators, and in particular \emph{bubblelators}, also known as gas evolution oscillators (cf.~\cite{SmithNoyes1983,BowersNoyes1983,Yuan1985,BarEli1992}). 
First, we suppose monomers are injected into the system at a constant source rate $S>0$:
\begin{equation}\label{e:Inj:ChemReact}
  \emptyset \xrightarrow{S} \{1\} \,.
\end{equation}
Second, we suppose clusters are removed at a rate proportional to a power law with removal coefficient $R>0$ and exponent $r\geq 0$ :
\begin{equation}\label{e:Depl:ChemReact}
  \{k\} \xrightarrow{R k^\rexp} \emptyset , \qquad k=2,3,\dots \, . 
\end{equation}
The resulting chemical reaction network~\eqref{e:BD:ChemReact}, \eqref{e:Inj:ChemReact}, \eqref{e:Depl:ChemReact} is open and no longer satisfies a detailed balance condition, in contrast to~\eqref{e:BD:ChemReact} alone.
By consequence, solutions may no longer dissipate free energy, and it becomes unclear whether long-time convergence to equilibrium always holds, or whether some more complicated dynamic behavior becomes possible. 
In this work we will provide evidence for the persistence of oscillations in time for a suitable approximate model of the network~\eqref{e:BD:ChemReact}, \eqref{e:Inj:ChemReact}, \eqref{e:Depl:ChemReact}.

\subsection{The classical Becker--D{\"o}ring model}

The Becker--D\"oring equations~\cite{BD1935} form an infinite system of kinetic equations that describes phase transitions in two-component mixtures where one of the phases has much smaller 
volume fraction than the other. In this case, the dilute phase consists of clusters of size $k \in \N$, where $k$ denotes the number of atoms, or monomers,
in the cluster. The main assumption in the Becker--D\"oring theory is that clusters evolve only by gain and loss of monomers. 
If $n_k$ denotes the density of clusters with $k$ monomers, and $J_k$ denotes the net rate of the reaction in \eqref{e:BD:ChemReact}, the equations read
\begin{align}
\partial_{t}n_{1} & =-J_{1}-\sum_{k=1}^{\infty}J_{k}\,, \label{bd1} \\
\partial_{t}n_{k} & =J_{k-1}-J_{k}\,, \qquad k\geq
2\,, \; \label{bd2} \\
J_{k}&=a_{k}n_{1}n_{k}-b_{k+1}n_{k+1}\,,\label{jkeq}
\end{align}
where $a_k, b_k$ are the respective attachment and detachment rate coefficients. 
The system of equations \eqref{bd1}-\eqref{jkeq} conserves the total mass $\rho$; that is,
\begin{equation}\label{masscons}
 \sum_{k=1}^{\infty} k n_k (t) = \sum_{k=1}^{\infty} kn_k(0)= \rho\,.
\end{equation}
Following work in statistical mechanics done by Penrose and collaborators \cite{Pen89,Pen97,PenEtal1984}
to model the dynamics of phase transitions, we take the coefficients to be of the form 
\begin{equation}\label{coeffassum}
a_{k}= k^{\alpha}\,, \qquad b_{k}=k^\alpha\Big(1+\frac{q}{k^{\gamma}}\Big ) \,, \qquad \text{ with }  q>0\,,\; \gamma \in (0,1)\,,\; \alpha \in [0,1)\,.
\end{equation}
The exponents $\alpha$ and $\gamma$ depend upon the geometry of clusters and the dominant mechanism of monomer transport:
For three-dimensional spheres dominated by diffusive transport,
 $\alpha=\sfrac{1}{3}$ and $\gamma=\sfrac{1}{3}$, while if cluster growth is limited by reactions
on the interface, we have $\alpha=\sfrac{2}{3}$ and $\gamma=\sfrac{1}{3}$. 
In the two-dimensional situation we have $\gamma=\sfrac{1}{2}$ and $\alpha=0$ and $\alpha=\sfrac{1}{2}$, respectively. 

The coefficient $q$ arises from the Gibbs-Thomson law and is proportional to surface tension. 
It plays a key role in determining a critical cluster size $k_{\rm crit}$ during the process of nucleation, 
a process which will prove fundamental throughout this paper.

We have chosen units for convenience such that the density of monomers 
in equilibrium with a planar phase boundary is 
\[
z_s = \lim_{k\to\infty}\frac{b_k}{a_k} = 1.
\]
This is also the maximum monomer concentration for which finite-mass steady states exist.  
The equilibrium state with this critical monomer concentration gives rise to a critical mass $\rho_s >0$ such that for any $\rho \in [0,\rho_s]$ an equilibrium solution $\bar n$ exists, which then has $\bar n_1\leq z_s$. 
For initial
data with supercritical mass $\rho>\rho_s$, it has been established in \cite{BCP86,Pen97} that the solution converges weakly in the long-time limit to
the equilibrium
solution with density $\rho_s$. The excess mass $\rho-\rho_s$ is transferred to larger and larger clusters as time proceeds and their evolution can be approximated by the classical LSW model for 
coarsening (see also \cite{Nie03,Schl19}). Furthermore, it has been shown in \cite{Pen89} that for certain initial data with small excess density it takes at least exponentially long time (in terms of the excess density)
until large clusters are created. Even though the proof is for specific data only, one expects that such metastable behavior appears for all generic data. (See also \cite{CDW95} for numerical
simulations.)

\subsection{The Becker--D{\"o}ring model with injection and depletion}\label{Ss.BD_inject_deplete}

In this paper we are interested in the Becker--D\"oring equations with injection of monomers and depletion of large clusters. 
More precisely for a given source rate $S>0$, removal rate $R>0$ and removal exponent $b\ge0$, we consider the system
\begin{align}
\partial_{t}n_{1} & =-J_{1}-\sum_{k=1}^{\infty}J_{k}+S\label{n1eq}\,,\\
\partial_{t}n_{k} & =J_{k-1}-J_{k}-R k^\rexp n_{k}\,, \qquad k\geq
2\,, 
\label{nkeq}
\end{align}
with $J_k$ as in \eqref{jkeq} and coefficients as in \eqref{coeffassum}.

It is well documented in the chemistry literature that temporal oscillations can persist in chemical--physical systems 
in which a phase transition creates strong nucleation peaks that lead to rapid growth of supercritical agglomerations which are later removed.
Specifically, the system~\eqref{n1eq}--\eqref{nkeq} has many similarities to models of {bubblelator} dynamics describing oscillatory
release of a gas (cf.~\cite{Morgan1916,BadgerDryden1939,SmithNoyes1983,PratsinisFriedlanderPearlstein1986,Friedlander2000}). For more background, we refer to Section~\ref{Ss.literature}.

Our goal is to obtain oscillatory solutions for an approximation of the model \eqref{n1eq}-\eqref{nkeq} under suitable
choices of the source term $S$ and the removal term $R$. 
The rough heuristics explaining the appearance of oscillatory solutions are the following. 
For small $S$ and $R$, the solution of \eqref{n1eq}, \eqref{nkeq} first evolves as in the classical model without source and removal terms---indeed, it evolves to 
a metastable steady state with locally constant nonzero fluxes $J_k$.
Slowly, the source term $S$ kicks in and raises the monomer density to a supercritical value with a small positive excess $n_1-1>0$ of order $\eps$ for small $\eps>0$.
This triggers the creation (or nucleation) of clusters larger than a critical size.
The nucleated supercritical clusters then grow by a process that depletes the monomer density
by a smaller amount of order $\eps^{\sfrac{1}{\gamma}}$ which is nevertheless enough to shut down large-cluster creation.
Large clusters eventually get destroyed, nucleation resumes, and the scenario repeats.

In order for this picture to be realized, the source term $S$ and the removal term $R$ have to be chosen in such a way that all the relevant terms balance. 
At steady state, source-driven nucleation balances removal of very large clusters. 
But the time it takes for nucleated clusters to grow large enough for effective removal 
introduces a delay that allows creation and removal to get out of phase.

\subsection{Limit model}\label{Ss.LimitModel}

In Section \ref{S.derivation} we identify the relevant scales, determine suitable $S$ and $R$ and derive a simplified model
formally valid in the limit $\eps\to0$.
This consists of an evolution equation for a rescaled monomer density $u$, defined in terms of the excess monomer density via 
\begin{equation}\label{eq:def:u}
n_1-1=\eps + \left(\frac{\eps}q\right)^{1/\gamma} u\,,
\end{equation}
and a transport equation for a rescaled density $f$ of large clusters with rescaled continuous size $x\in (0,\infty)$.
Section~\ref{Ss.scales} contains the precise definition of the scales and a rescaled removal rate $\rfac>0$.
Taking those for granted, we formally derive the system
\begin{align}
	\partial_\tau u(\tau) &= 1 - \int_0^\infty x^\alpha f(x,t) ,dx \,, \label{eq:limit:u} \\
	\partial_{\tau}f(x,\tau)  +\partial_{x}\big(  x^{\alpha}f(x,\tau)  \big)   &  =-\rfac  
	x^\rexp f(x,\tau)  \,, \qquad  x>0\,, \label{eq:limit:transport} \\ 
	x^\alpha f(x,\tau) &\to e^{u(\tau)} \, , \qquad x\to 0 \,.\label{eq:limit:boundary} 
\end{align}
The key ingredients that go into the derivation of this system are the following. First, the distribution of clusters
for $k$ around the critical size $k_{\crit}$ or smaller is taken as quasistationary, corresponding to constant-flux
states for the Becker-D\"oring equations, which are parametrized by the monomer density.  
The transport equation arises by a continuum approximation to the difference equation \eqref{bd2}
in the supercritical range $x\simeq (k-k_{\rm crit})/X >0$, with size scale $X$ exponentially large in $\eps$.
The two size regimes are related by matching fluxes around $k\simeq k_{\rm crit}$.
For the quasistationary states, a continuum approximation results in an Arrhenius law 
giving the exponential dependence on $u$ in \eqref{eq:limit:boundary}. 
The precise scaling of $u$ by $\eps^{\sfrac{1}{\gamma}}$ leads to 
a change of order one in the Arrhenius factor in the boundary condition~\eqref{eq:limit:boundary} 
when $u$ has a change of order one. 

We remark that a simple computation yields the mass balance law
\begin{equation}\label{massconservation}
	\partial_{\tau}\left(  u(\tau)+\int_{0}^{\infty}x f(x,\tau) \, dx\right)
	=1-\rfac \int_{0}^{\infty}x^{\rexp+1} f(x,\tau) \, dx .
\end{equation}
Hence, the total mass in the system increases though the influx of monomers in~\eqref{eq:limit:u} at unit rate,
and decreases due to the removal of large clusters on the right hand side of~\eqref{eq:limit:transport}.

\subsection{Oscillations via Hopf bifurcation}

For the approximate model~\eqref{eq:limit:u}--\eqref{eq:limit:boundary}, time-periodic solutions satisfy a delay-differential equation with infinite delay horizon.
To our knowledge, a mathematically rigorous Hopf bifurcation theorem has been proven only for finite delay horizons.
Nevertheless, we argue in Section~\ref{S.hopf} that Hopf bifurcations from stationary solutions should occur as one varies the removal parameter $\rfac$.
In particular, we make a careful analysis of the spectrum of the linearized problem around a steady state
solution and track the dependence of the eigenvalues of the linearized operator upon the removal parameter $\rfac$. 
We identify points where eigenvalues cross the imaginary axis transversely,
and we perform a formal expansion to determine the direction that bifurcation should occur, 
which should indicate when \emph{stable} periodic solutions appear. 
We provide analytical and numerical evidence regarding bifurcation points and their dependence on $\alpha$ and $\rexp$ in Section~\ref{s:parms}.

\subsection{Oscillations via desycronization of source and removal}\label{Ss.description}

In Figure~\ref{fig:numeric} a numerical solution of the limit model \eqref{eq:limit:u}--\eqref{eq:limit:boundary} is shown. 
The oscillations are large amplitude, indicating that the chosen parameters are already far beyond the Hopf bifurcation point. 
The flux profiles in this regime develop an interesting structure involving the transport of rather sharp peaks,
which one can understand in a physical way that we wish to explain. 
\newcommand{\subfigheight}{0.25\textwidth}
\begin{figure}[!ht]
	\centering
	\begin{tikzpicture}
%		\begin{axis}[colormap/viridis,width=0.5\textwidth,height=0.33\textwidth, grid=major, axis y line=left, axis x line=top,xticklabel style={overlay},
		\begin{axis}[colormap/blackwhite, %viridis
			width=0.5\textwidth,height=\subfigheight, grid=major, axis y line=left, axis x line=top,xticklabel style={overlay},
			xmin=0,xmax=600, ymin=-22, ymax=0, xlabel={$\tau$}, %ylabel={$u$}, , xlabel={$\tau$}
			every axis x label/.style={at={(current axis.right of origin)},anchor=west}]
			\addplot+[thick, black, line join=round, mark=none] table [x=t,y=u]{simu-u.csv};\label{p:u}
			\addplot+[only marks, mark=*, scatter, point meta=x, point meta min=523, point meta max=563] table [x=t,y=u]{simu-u-marks.csv};\label{p:u-marks}
		\end{axis}
	\end{tikzpicture}\hfill%
	\begin{tikzpicture}
		%semilogyaxis
		\begin{semilogyaxis}[colormap/viridis,cycle list={[colors of colormap={0,250,500,750,1000}]},
			width=0.5\textwidth, height=\subfigheight, grid=major, axis lines=left, xmin=0,xmax=1000, ymax=1, ymin=5e-13,xticklabel style={overlay}, %ylabel={$h$},, 
			xlabel={$x$},
			log origin=infty,every axis x label/.style={at={(current axis.right of origin)},anchor=west}]
			\addplot+[thick] table [x=x,y=t523]{simu-f.csv};\label{p:t523}
%			\addplot+[thick] table [x=x,y=t533]{simu-f.csv};\label{p:t533}
%			\addplot+[thick] table [x=x,y=t543]{simu-f.csv};\label{p:t543}
%			\addplot+[thick] table [x=x,y=t553]{simu-f.csv};\label{p:t553}
%			\addplot+[thick] table [x=x,y=t563]{simu-f.csv};\label{p:t563}
		\end{semilogyaxis}
	\end{tikzpicture}\hspace{1em}

\bigskip
\begin{tikzpicture}
	\begin{semilogyaxis}[colormap/viridis,cycle list={[colors of colormap={0,250,500,750,1000}]},
		width=0.5\textwidth, height=\subfigheight, grid=major, axis lines=left, xmin=0,xmax=1000, ymax=1, ymin=5e-13, xticklabel style={overlay},%ylabel={$h$},
		log origin=infty,every axis x label/.style={at={(current axis.right of origin)},anchor=west}]
%		\addplot+[thick] table [x=x,y=t523]{simu-f.csv};\label{p:t523}
		\addplot+[thick] table [x=x,y=t533]{simu-f.csv};\label{p:t533}
%		\addplot+[thick] table [x=x,y=t543]{simu-f.csv};\label{p:t543}
%		\addplot+[thick] table [x=x,y=t553]{simu-f.csv};\label{p:t553}
%		\addplot+[thick] table [x=x,y=t563]{simu-f.csv};\label{p:t563}
	\end{semilogyaxis}
	\end{tikzpicture}\hfill
\begin{tikzpicture}
	\begin{semilogyaxis}[colormap/viridis,cycle list={[colors of colormap={0,500,750,1000}]},
	width=0.5\textwidth, height=\subfigheight, grid=major, axis lines=left, xmin=0,xmax=1000, ymax=1, ymin=5e-13, xticklabel style={overlay},%ylabel={$h$},
	log origin=infty,every axis x label/.style={at={(current axis.right of origin)},anchor=west}]
%		\addplot+[thick] table [x=x,y=t523]{simu-f.csv};\label{p:t523}
%		\addplot+[thick] table [x=x,y=t533]{simu-f.csv};\label{p:t533}
		\addplot+[thick] table [x=x,y=t543]{simu-f.csv};\label{p:t543}
%		\addplot+[thick] table [x=x,y=t553]{simu-f.csv};\label{p:t553}
%		\addplot+[thick] table [x=x,y=t563]{simu-f.csv};\label{p:t563}
\end{semilogyaxis}
\end{tikzpicture}\hspace{1em}

\bigskip

\begin{tikzpicture}
\begin{semilogyaxis}[colormap/viridis,cycle list={[colors of colormap={0,750,1000}]},
	width=0.5\textwidth, height=\subfigheight, grid=major, axis lines=left, xmin=0,xmax=1000, ymax=1, ymin=5e-13, xticklabel style={overlay}, %ylabel={$h$},
	log origin=infty,every axis x label/.style={at={(current axis.right of origin)},anchor=west}]
%		\addplot+[thick] table [x=x,y=t523]{simu-f.csv};\label{p:t523}
%		\addplot+[thick] table [x=x,y=t533]{simu-f.csv};\label{p:t533}
%		\addplot+[thick] table [x=x,y=t543]{simu-f.csv};\label{p:t543}
		\addplot+[thick] table [x=x,y=t553]{simu-f.csv};\label{p:t553}
%		\addplot+[thick] table [x=x,y=t563]{simu-f.csv};\label{p:t563}
\end{semilogyaxis}
\end{tikzpicture}\hfill
\begin{tikzpicture}
\begin{semilogyaxis}[colormap/viridis,cycle list={[colors of colormap={0,1000}]},
	width=0.5\textwidth, height=\subfigheight, grid=major, axis lines=left, xmin=0,xmax=1000, ymax=1, ymin=5e-13, xticklabel style={overlay}, %ylabel={$h$},
	log origin=infty,every axis x label/.style={at={(current axis.right of origin)},anchor=west}]
%		\addplot+[thick] table [x=x,y=t523]{simu-f.csv};\label{p:t523}
%		\addplot+[thick] table [x=x,y=t533]{simu-f.csv};\label{p:t533}
%		\addplot+[thick] table [x=x,y=t543]{simu-f.csv};\label{p:t543}
%		\addplot+[thick] table [x=x,y=t553]{simu-f.csv};\label{p:t553}
		\addplot+[thick] table [x=x,y=t563]{simu-f.csv};\label{p:t563}
\end{semilogyaxis}
\end{tikzpicture}\hspace{1em}
\medskip

	\caption{Time evolution of $u$ (top left), and semi-log plots of flux $x^\alpha f$ 
at times $\tau=523$~(top right), $\tau=533$ and $543$ (middle row), $\tau=553$ and $563$ (bottom row). The times $\tau=523,\dots,563$ are indicated by circle marks in the top left plot.   %\newline
		Parameters: $\rfac=0.1$, $\alpha=\frac{1}{3}$, $\gamma=\frac{1}{3}$, $b=\frac{2}{3}$.}\label{fig:numeric}
\end{figure}
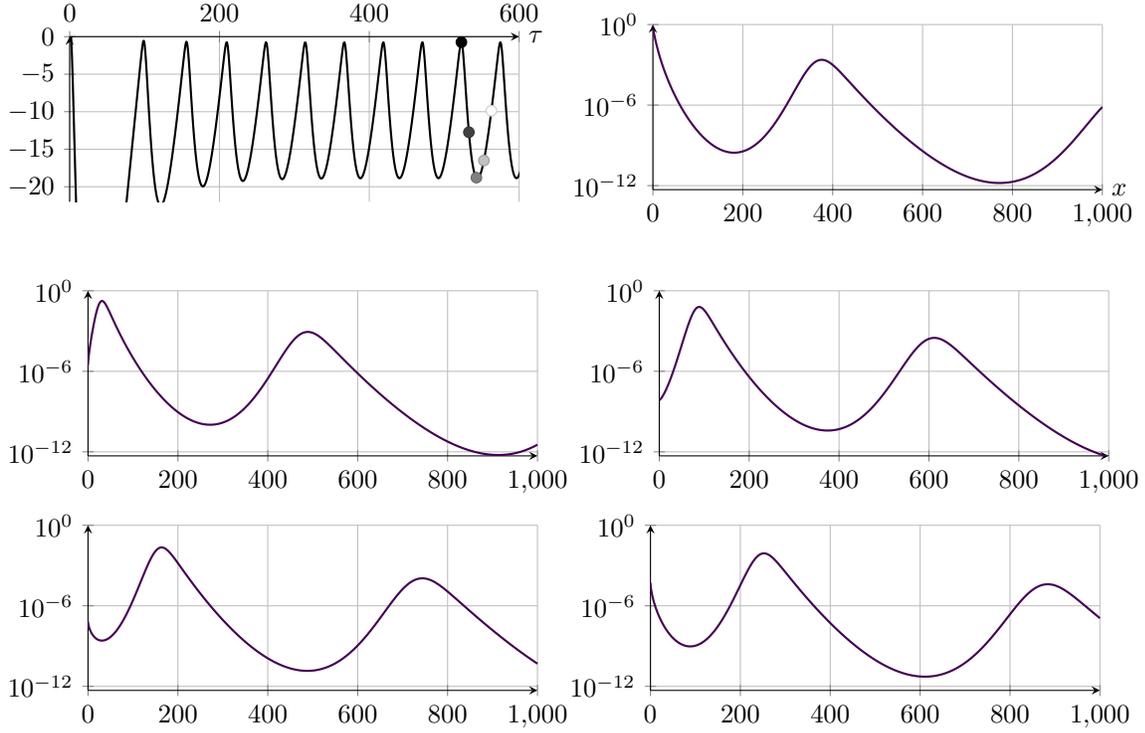

Before discussing that, we remark that stable oscillations should also exist in the full Becker--D{\"o}ring model with injection and depletion~\eqref{n1eq}, \eqref{nkeq}. 
The numerical computation of such oscillations seems to be a challenging matter, however, due to the multiscale nature of the Becker--Döring system, as already observed in~\cite{CDW95}. 
The main difficulty is that the scales associated to the problem depend exponentially on the small parameter $\eps$. For instance, we obtain for $\eps=0.1$ that the typical size of clusters that are involved in the dynamics of the system~\eqref{n1eq}, \eqref{nkeq} is of order $10^{12}$ with a typical time-scale of order $10^9$ (cf. Subsection~\ref{Ss.scales}).
For larger values of $\eps$ a numerical approach might be feasible, but the formal asymptotic approximation done in this paper may not apply.

\begin{figure}[!ht]
 \centering
 \setlength{\unitlength}{0.75\textwidth}%
  \begin{picture}(1,0.56020297)%
    \put(0,0){\includegraphics[width=\unitlength]{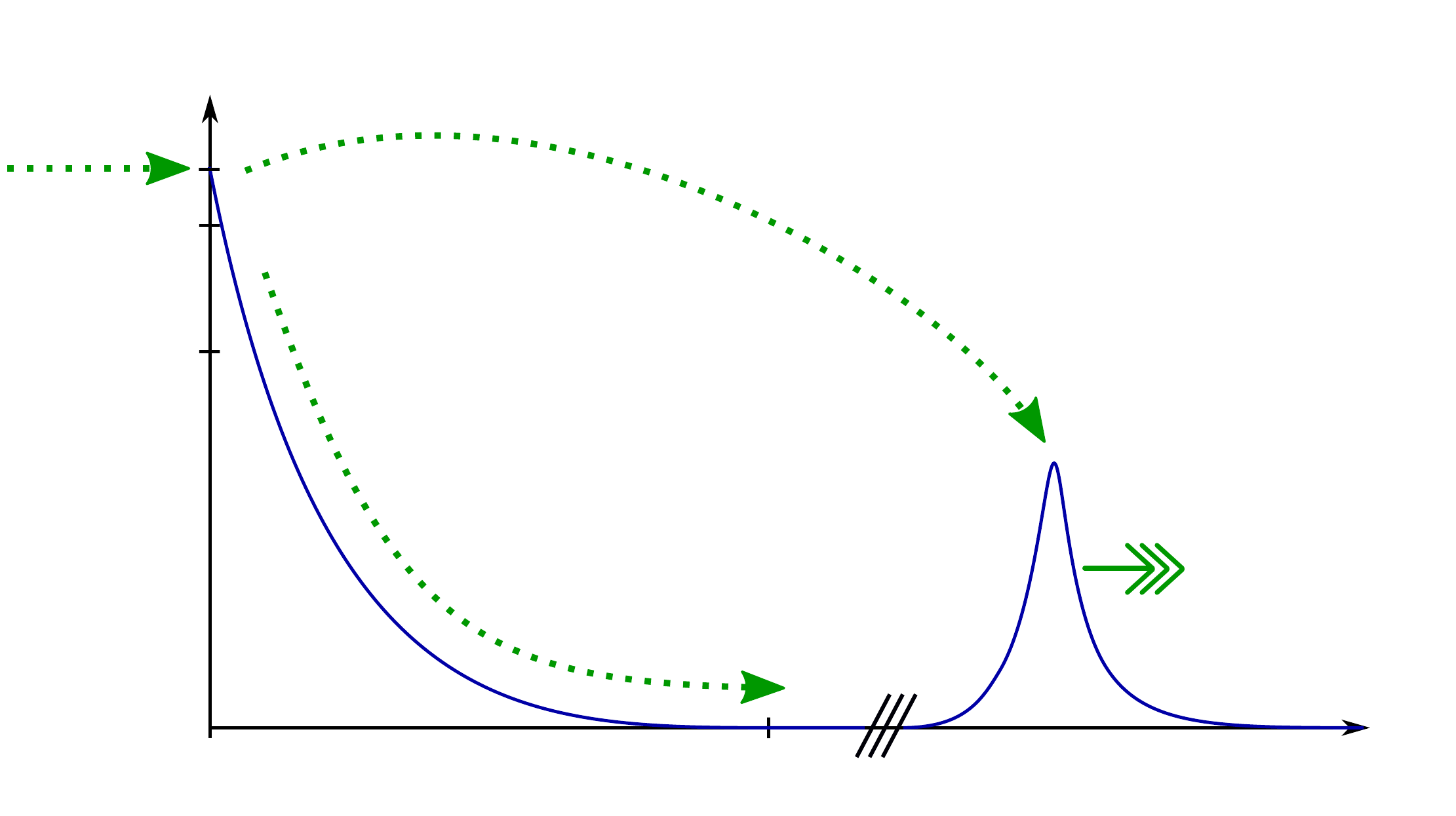}}%
    \put(0.11166376,0.30779252){$1$}%
    \put(0.94262686,0.05074474){$k$}%
    \put(0.13,0.5){$n_k$}%
    \put(0.05,0.4){$1+\eps$}%
    \put(0.13004492,0.02250351){$1$}%
    \put(0.51263459,0.02350007){$k_{\crit}$}%
    \put(0.7,0.02350007){$X$}%
    \put(0.14906643,0.44709719){$u$}%
    \put(0.2664336,0.19426132){$J_{\nucl}$}%
    \put(0.74068943,0.1977548){$V_{\transport}$}%
    \put(0.4119297,0.46142681){$J_{\cond}$}%
    \put(0.05,0.45){$S$}%
  \end{picture}%
 \caption{Illustration of transport mechanisms (not to scale). For explanation, see text.}\label{fig:flux}
\end{figure}
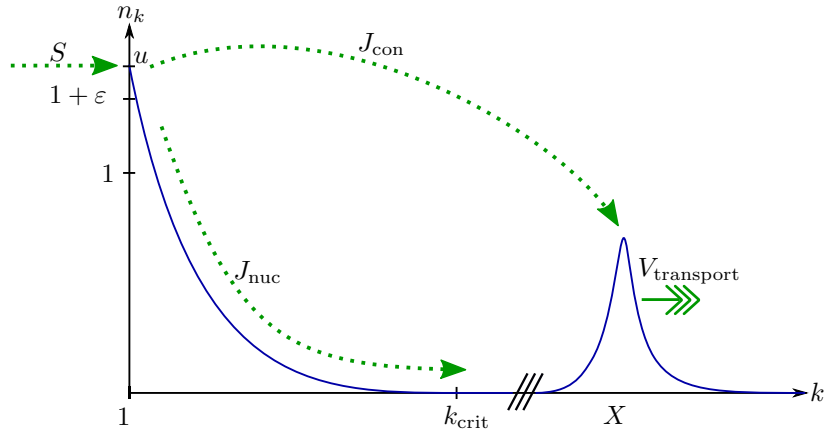

In Figure~\ref{fig:flux}, different fluxes in the Becker--D{\"o}ring model relevant for the description of the oscillation mechanism are shown, superimposed on a schematic and highly exaggerated plot of $n_k$ {\it vs} $k$. The
first crucial quantity is the \emph{critical size} $k_{\crit}$, as in~\cite{Pen89}. This depends on the monomer excess $\eps = n_1{-}1$, and is defined for simplicitiy here as the value of $k$ for which $a_k n_1 - b_k$ vanishes (cf.~Section~\ref{Ss.steadystates}). For the rates~\eqref{coeffassum} and recalling that $z_s=1$, it holds 
\begin{equation}\label{criticalsize}
 k_{\crit} = \Big( \frac{q}{\eps}\Big)^{\frac{1}{\gamma}} \,.
\end{equation}
In the following, the critical size $k_\crit$ is used to distinguish small from supercritical clusters. It also provides the relevant next scale for the monomer expansion in~\eqref{eq:def:u}, since $n_1 -1 = \eps + u/k_\crit$.
The first driving mechanism is the flux of mass through the small clusters to beyond the critical size $k_{\crit}$, denoted by $J_{\nucl}$. In the literature (cf.~Friedlander~\cite{Friedlander2000}), this process is called \emph{homogeneous nucleation}. 
Because this process is diffusion-dominated for $k$ around $k_{\crit}$, we obtain a boundary layer with size of order $\eps^{-\sfrac{1}{\gamma}}$, 
resulting in an Arrhenius relation for the flux. This is reflected in the limit model through the exponential boundary condition~\eqref{eq:limit:boundary}. 

A second flux $J_{\cond}$ depletes the monomer concentration through the mechanism of large-cluster growth by direct absorption of monomers, 
also called \emph{condensation}. This flux is reflected in the limit model as the integral loss term in \eqref{eq:limit:u}, and also leads to the transport term in~\eqref{eq:limit:transport}.

These two fluxes need to be balanced with the source term $S$, and the time-scale chosen accordingly, which is done in a careful analysis of the scales in Section~\ref{Ss.scales}. 

Both systems, the Becker--Döring model~\eqref{n1eq}--\eqref{nkeq} and the limit model~\eqref{eq:limit:u}--\eqref{eq:limit:boundary}, allow for stationary states with time-independent concentration of monomers, and we believe that these are stable for large depletion parameters $R$ or $\rfac$, respectively. 
However, our results indicate that time-periodic solutions of the limit model exist for small particle-removal rates. 

An explanation is that in this case clusters must grow very large for effective removal.
This depletes the monomer concentration due to the condensation flux $J_{\cond}$ that drives the growth of large clusters,
and leads to time delays in replenishing those clusters. 
It takes more time both to restart nucleation from a lower monomer level 
and to grow supercritical clusters to sizes large enough for removal.
In this way the generation and removal of large clusters can \emph{desynchronize}, resulting in temporal oscillations.

In somewhat more detail, the mechanism for oscillations works as follows.
Peaks in the size distribution of large clusters form through a process mediated by the 
sensitive (exponential) dependence of the nucleation rate upon the concentration of monomers:
A sufficient excess of monomers above the critical concentration $n_1=1$ triggers rapid growth of the number of supercritical clusters through nucleation.
The creation of enough supercritical clusters then produces a large condensation flux $J_{\cond}$ 
which forces the concentration of monomers to decrease in spite of the source term.
This stops, or drastically slows, the nucleation reaction transferring mass through the critical size.

The peak of supercritical clusters generated in this way then continues to be 
transported to ever-larger sizes by the condensation mechanism, which continues to consume monomers.
At sufficiently large sizes, the rate that clusters are removed from the system becomes dominant and the peak is eliminated.
This makes the condensation flux $J_{\cond}$ small again, which allows the source term to force the monomer concentration higher, and start the cycle over.

\subsection{Sharp peak model} 

The mechanism for oscillations just described can be implemented in a simple model with only two elements, 
namely the excess monomer density $u(t)$ and the size distribution $f(x,t)$ for clusters of supercritical size $x>x_{\crit}$. 
The clusters increase their size at some constant rate by absorbing monomers, 
and they are removed instantaneously when size reaches some terminal value $x_{\rem}>x_{\crit}$.
If the source of monomers drives their excess to reach a certain nucleation threshold $u_{\nucl}$, 
a sudden \emph{sharp peak} of large clusters is created just above the critical size.

The balance law for the monomer excess takes the form
\begin{equation}\label{peak:u}
\partial_{t}u(t)=1 - \int_{x_{\crit}}^{x_{\rem}} f(x,t)\, dx ,
\end{equation}
and one starts with initial excess $u(0)<u_{\nucl}$ below the nucleation threshold.
The size distribution~$f$ is advected at constant speed $1$, satisfying the transport equation
\begin{equation}\label{peak:transport}
\partial_{t}f+\partial_{x}f=0 ,\qquad x_{\crit}<x<x_{\rem} .
\end{equation}
We specify a zero influx condition $f(x_{\crit},t)=0$ as long as the monomer flux remains below threshold, i.e., $u(t)<u_{\nucl}$.
When the monomer excess reaches the threshold at some time $t_*$, however, 
a delta-mass concentration of supercritical clusters is instantly nucleated at $x=x_{\crit}$, giving the jump condition
\begin{equation}\label{peak:boundary}
 f(\cdot,t_*^+) = f(\cdot,t_*^-) + f_0 \delta_{x_{\crit}}  \qquad\mbox{if} \quad u(t_*)=u_{\nucl}.
\end{equation}
For $f_0>1$ and $f\equiv 0$ initially, say, this model produces a sawtooth evolution for $u$,
with $\partial_t u = 1>0$ during time intervals when no supercritical clusters exist in the system,
and $\partial_t u=1-f_0<0$ during intervals after a peak of clusters has been nucleated and before it is removed
upon reaching the outflow boundary $x=x_{\rem}$.

The resulting oscillatory evolution in this model clearly illustrates
desynchronization of supercritical cluster generation and removal.
The first two equations~\eqref{peak:u} and~\eqref{peak:transport} of this model are very similar to~\eqref{eq:limit:u} and~\eqref{eq:limit:transport} by setting $\alpha=0$ and letting $b\to \infty$.
The only main difference is  that the exponential boundary condition~\eqref{eq:limit:boundary} is changed to the jump condition~\eqref{peak:boundary}, leading here to the periodic production of sharp peaks.

\subsection{Related literature}\label{Ss.literature}

Oscillations in chemical reaction networks have been well known ever since the Belousov--Zhabotinsky reaction was described~\cite{Belousov1959,Zhabotinsky1964} and the \emph{Brusselator} found by Prigogine and Lefever~\cite{PrigogineLefever1968} (see~\cite{Tyson1976} for the name). Later, the \emph{Oregonator} was introduced by Field and Noyes~\cite{FieldNoyes1974} as a simpler model that develops temporal oscillations involving only five species.
The mathematical analysis of these systems reveals that the basic mechanism behind the oscillations is of Lotka--Volterra type~\cite{Lotka1910,Lotka1920,Volterra1926}.

In contrast to this, the model from Section~\ref{Ss.BD_inject_deplete} is motivated by mechanisms
found in the dynamics of \emph{bubblelators}, also known as gas evolution oscillators~\cite{SmithNoyes1983,BowersNoyes1983,Yuan1985,BarEli1992}. 
These are chemical--physical systems~\cite{SmithNoyesBowers1983} in which, due to some reaction mechanism, a dissolved gas is constantly generated in a solvent, 
leading to a steady increase in supersaturation and an eventual burst of nucleation and growth of gas bubbles.
The first experimental report of such a system is ascribed to Morgan~\cite{Morgan1916}, who observed an oscillatory release of gas during
dehydration of formic acid in concentrated sulfuric acid. 
The dynamics and growth of bubbles resembles a mechanism similar to the Becker--D{\"o}ring paradigm: 
The growth of bubbles is effected mainly through the absorption of gas emerging from the supersaturated solution into the expanding bubbles. 
Upon discretization in size, this suggests that the physical growth mechanism is reflected by fluxes resembling those in~\eqref{jkeq}. 
Lastly, large bubbles randomly leave the solvent, depending on the setup of the experiment, mainly due to buoyancy, 
and this is loosely reflected in the system~\eqref{nkeq} by the removal term $-Rk\rexp n_k$. 

More experimental evidence of oscillatory concentrations in a mixture of nitric
oxide and coal gas was obtained by Badger and Dryden~\cite{BadgerDryden1939}. A theoretical model of
Becker--D\"{o}ring type was proposed by Pratsinis, Friedlander, and Pearlstein~\cite{PratsinisFriedlanderPearlstein1986} (cf.\ also \cite{Friedlander2000}). The
model in \cite{PratsinisFriedlanderPearlstein1986} has some physical similarities with the limit model
obtained in Section~\ref{Ss.LimitModel}. More precisely, the model uses the
exponent $\alpha=\sfrac{2}{3}$ in~\eqref{coeffassum} and $b=0$ in~\eqref{nkeq}.
In this case, it is possible to obtain a closed system of ODEs for
the evolution of the three lowest-order moments of the distribution of radii of supercritical clusters.
This system of ODEs is coupled with the concentration of monomers $n_{1}$ by means of an Arrhenius formula yielding
the nucleation rate of supercritical clusters as a function of the monomer concentrations. 
The resulting model in \cite{PratsinisFriedlanderPearlstein1986} is a system of ODEs for which the existence of periodic solutions 
is demonstrated using a Hopf bifurcation argument. 
A similar reduction of a Becker--D\"{o}ring model to a system of ODEs which has periodic solutions can be found in \cite{McGrawSaunders1984}. 
We can show (see the Appendix~\ref{Ss.momentmodels})
that, for some particular choices of exponents, 
the limit model from Section~\ref{Ss.LimitModel} can be reduced to a system of ODEs 
which has the same structure as those obtained in~\cite[p.~293ff]{Friedlander2000}.

A boundary condition similar to~\eqref{eq:limit:boundary} was derived by Farjoun and Neu~\cite{FarjounNeu2008}
in a physical study describing the depletion of a supercritical concentration of monomers (without source) 
due to the nucleation of supercritical clusters. 
The nucleation rate is approximated using an Arrhenius law (or Frank-Kamanetskii) type of formula. 
A key observation made in \cite{FarjounNeu2008} in the derivation of this
boundary condition is that small changes in the concentration of monomers
$n_{1}$ yield significant changes in the nucleation rate.
This same point underpins the present study.

Recently, oscillations for a Becker--D{\"o}ring model with atomization were proved to exist by two of the present authors in \cite{PV19}. 
The model in~\cite{PV19} is closed and has no external source or removal terms ($S=r=0$ in~\eqref{n1eq}--\eqref{nkeq}).
The atomization of clusters of a maximal size $M$ into $M$ monomers provides a closed feedback mechanism from large clusters to monomers,
which could be considered to replace injection and depletion.
This model has a Hopf bifurcation for suitable small atomization rate, when $M$ is large.

In the physical literature, temporal oscillations in coagulation-fragmentation models permitting interactions of clusters of any sizes
(thus not of Becker-D\"oring type) have been reported in numerical simulations, 
by R.~C.~Ball et al.~\cite{Ball2012} for cases with monomer injection and cluster removal above a fixed size,
and in the works \cite{Matveev2017,Brilliantov2018} for cases incorporating a nonlinear atomization mechanism.
The onset of a Hopf bifurcation for a model consisting of coagulation with monomers and atomization has been recently shown numerically in~\cite{BudzinskiyMatveevKrapivsky2021}.

Lastly, Doumic et al.~\cite{Fellner_etal2018} consider a model for prion dynamics of Becker--D{\"o}ring type, which exhibits very slowly damped oscillations. 
The model in~\cite{Fellner_etal2018} assumes that the polymer chains interact with two types of chemicals yielding respectively 
increase and decrease of the length of the polymer chain. These chemicals interact between themselves by means of a modified Lotka--Volterra type of equation, 
which is coupled with the Becker--D{\"o}ring part of the system. It is well known that Lotka--Volterra models may admit periodic solutions. 
The interaction of the Lotka--Volterra equation with the Becker--D{\"o}ring part is responsible for the damping of the oscillations observed in~\cite{Fellner_etal2018},
but the specific form in which this damping takes place is not well understood at the moment.

\section{Derivation of the limit model}\label{S.derivation}

In this section, we describe a regime in which the discrete Becker--Döring model with source  and depletion~\eqref{n1eq}--\eqref{nkeq} can be formally approximated by the limit model~\eqref{eq:limit:u}--\eqref{eq:limit:boundary}. 

Our analysis is based on understanding the interaction between small and large clusters, where the separation scale is given by the critical size~\eqref{criticalsize}. 
We shall argue that for clusters smaller or comparable to $k_\crit$, solutions are close to a constant-flux steady state of the Becker--Döring equation parametrized by the monomer concentration, and we derive an Arrhenius law relating the flux and monomer concentration.
We provide additional discussion of the quasistationary assumption in 
\AppendixSupp~\ref{sec:quasistationary},
where we describe the time scale for a drift-diffusion approximation in the boundary layer where $k$ is near critical size.

In Section~\ref{Ss.transport}, we show that the dynamics of large clusters are well approximated by a transport equation, which will rescale to~\eqref{eq:limit:transport}. 
The boundary condition for this transport equation is obtained in Section~\ref{Ss.fluxmatch} by matching with the fluxes from Section~\ref{Ss.steadystates}.
Then in Section~\ref{Ss.monomers}, we study the evolution of the monomer fluctuation $u$ in~\eqref{eq:def:u}. This is a balance between (i) the nucleation flux towards large clusters through the critical size $k_\crit$ and (ii) a condensation flux due to the growth of large clusters. 

Finally, in Section~\ref{Ss.scales}, we will identify a time scale $T$, a macroscopic cluster size scale $X$ much larger the critical cluster size~\eqref{criticalsize}, a macroscopic cluster density $F$, as well as the source rate $S$ and the depletion rate $R$ for which the Becker--Döring system~\eqref{n1eq}--\eqref{nkeq} with injection and depletion can be approximated by the limit model~\eqref{eq:limit:u}--\eqref{eq:limit:boundary}.

\subsection{Steady states with constant flux}\label{Ss.steadystates}

We recall results from \cite{Pen89} about steady states with constant flux.
We consider in particular the case of small excess density and find suitable asymptotic expressions of the steady states in this regime.
It is convenient to recall first the formulas for equilibrium solutions with zero flux, given by 
\begin{equation}\label{equil1}
 {\bar n}_k= Q_k {\bar n}_1^k\,, \qquad 
 Q_k=\prod_{l=1}^{k-1}\frac{a_l}{b_{l+1}} \,,
\end{equation}
where the parameter ${\bar n}_1$ represents the equilibrium concentration of monomers. With coefficients as in \eqref{coeffassum} we find 
\begin{equation}\label{qasymp}
 \log\bigl( a_k Q_k \bigr) = \log \biggl( \prod_{l=1}^k \frac{a_l}{b_l} \biggr) 
 = \sum_{l=1}^k \log\biggl(\frac{1}{1+\frac{q}{l^{\gamma}}}\biggr) 
 \simeq -\frac{qk^{1-\gamma}}{1-\gamma}(1+o(1))
\qquad\text{as $k\to\infty$.}
\end{equation}
Thus we see that the series
$\sum_{k=1}^{\infty} k Q_k {\bar n}_1^k$, which represents the mass of $({\bar n}_k)$,  
has radius of convergence $1$ and converges for ${\bar n}_1=1$. 
The corresponding critical mass is denoted as $\rho_s =\sum_{k=1}^{\infty} kQ_k$.

For super-critical monomer density  ${\bar n}_1>1$, the zero-flux equilibrium solution $Q_k{\bar n}_1^k$ grows exponentially at infinity and those solutions will not play a role. 
As in \cite{Pen89} we consider in this regime steady states  with constant flux, where the flux is chosen such that the steady state remains bounded as $k\to \infty$. 
More precisely, we look for given $\bar n_1$ for bounded solutions $\{N_k(\bar n_1)\}_{k=\N}$ of 
\begin{equation}\label{Jstateq}
N_1= \bar n_1 \qquad \text{and}\qquad a_{k{-}1}{\bar n}_1 N_{k{-}1}({\bar n}_1) - b_k N_k({\bar n}_1)= J({\bar n}_1)\,, \qquad k\geq2\,,
\end{equation}
where $J(\bar n_1)$ is part of the unknown. 
It has been shown in \cite[Lemma 1]{Pen89} that for each ${\bar n}_1>1$ there exists a unique solution to this problem, given by the formula
\begin{equation}\label{stst}
 N_k({\bar n}_1)= J({\bar n}_1) Q_k {\bar n}_1^k \sum_{l=k}^{\infty} \frac{1}{a_l Q_l {\bar n}_1^{l + 1}}
 \qquad\text{where}\qquad 
 \frac 1{J({\bar n}_1)}  ={\sum_{l=1}^{\infty}\frac{1}{a_{l}Q_{l} {\bar n}_1^{l+1}}}\,.
\end{equation}
Furthermore, for fixed ${\bar n}_1$,  $a_k N_k({\bar n}_1)$ decreases monotonically with $k$ (which implies that $N_k({\bar n}_1)$ is bounded and $N_k({\bar n}_1) \to 0$ as $k \to \infty$ if $\alpha>0$), while for fixed $k$, $\tfrac{N_k({\bar n}_1)}{{\bar n}_1}$ 
increases monotonically with ${\bar n}_1$  and we have 
for ${\bar n}_1>1$ that $Q_k < N_k({\bar n}_1) < Q_k {\bar n}_1^k$ for $k \in \N$.

\subsubsection*{Asymptotics of steady states and flux}
We are particularly interested in the asymptotic behavior of $\{N_k(\bar n_1)\}_{k=\N}$ and $J(\bar n_1)$ for slightly supercritical density  ${\bar n}_1>1$. Thus we introduce the small parameter
\begin{equation}\label{epsdef}
\eps = {\bar n}_1-1 \ll 1\, .
\end{equation}
We shall argue that the asymptotics for the constant flux are given by 
\begin{equation}\label{Jasymp}
	J(\bar n_1) \simeq
	J_\infty := 
	\sqrt{\frac{\gamma}{2\pi q^{\frac1\gamma}}}
	\eps^{\frac{\gamma+1}{2\gamma}}
	\exp\biggl(-\frac{\gamma}{1-\gamma} q^{\frac{1}{\gamma}} \eps^{-\frac{1-\gamma}{\gamma}}  \biggr) 
	\qquad\text{as } \eps\to 0,
\end{equation}
and for the corresponding steady states by
\begin{equation}\label{nstatasymptotics}
	N_k({\bar n}_1) \simeq \frac{J({\bar n}_1)}{\eps a_k} = \frac{J({\bar n}_1)}{\eps k^{\alpha}} \qquad \text{ for } k \gg k_{\crit} = \Big( \frac{q}{\eps}\Big)^{\frac{1}{\gamma}} \,.
\end{equation}

{\em Derivation of \eqref{Jasymp}.}
The critical cluster size $k_\crit$ from~\eqref{criticalsize} is a crucial quantity occurring in the analysis of this paper, as in \cite{Pen89}. 
Improving on \eqref{qasymp}, we can write
\begin{equation}\label{Gkdef}
\frac 1{a_kQ_k {\bar n}_1^k} = e^{G(k)}(1+o(1))\,,
\qquad \text{with} \quad 
G(k)= -k\log\bar n_1 + \int_1^k \log(1+q l^{-\gamma})\,d l + C \,.
\end{equation}
The series for $\sfrac{1}{J(\bar n_1)}$ in \eqref{stst}
is dominated by terms with $l$ near the point where $k\mapsto G(k)$ is maximized, and it happens exactly at 
$k_\crit=(\sfrac{q}{\eps})^{\sfrac1\gamma}$,
which is consistent with~\eqref{criticalsize}.
Laplace's method provides an approximation to the series: Noting that 
\[
G''(k) = \frac{-\gamma q k^{-\gamma-1}}{1+q k^{-\gamma}}
\]
and $G'''(k)/G''(k)=O(\sfrac1k)$, the expansion $ G(k) = G(k_\crit)+\frac12 G''(k_\crit)(k-k_\crit)^2(1+o(1)) $
is valid for $|k-k_\crit|<k_\crit^p$ for any $p<1$.  
Choosing $p>\frac12(1+\gamma)$ allows $k_\crit^{-\gamma-1}(k-k_\crit)^2$ to be large, whence we find
\[ 
\sum_{|l-k_\crit|<k_\crit^p}
\frac{1}{a_l Q_l {\bar n}_1^{l}}  \simeq 
e^{G(k_\crit)} 
\int_{-\infty}^\infty 
e^{ \frac12 G''(k_\crit)(k -k_\crit)^2}
\,dk
= e^{G(k_\crit)}\sqrt{\frac{2\pi}{-G''(k_\crit)}}
\qquad \text{as $\eps\to0$.  }
\]
The remaining part of the series for $\sfrac{1}{J(\bar n_1)}$ is small relative to this, 
and since $n_1\simeq 1$ it follows~\eqref{Jasymp} from $J(\bar n_1) \simeq \sqrt{\frac{-G''(k_\crit)}{2\pi}} e^{-G(k_\crit)}$ and by noting that 
\begin{equation}\label{Gmax}
G(k_\crit) \simeq \frac{qk_\crit^{1-\gamma} }{1-\gamma} (1+o(1))-k_\crit \log\bar n_1 \simeq
\frac{k_\crit\eps \gamma}{1-\gamma} =
 \frac{\gamma}{1-\gamma} q^{\frac1\gamma} \eps^{1-\frac1\gamma}.
\end{equation}

{\em Derivation of~\eqref{nstatasymptotics}.}
Similarly we obtain by recalling the relation~\eqref{Gkdef} for $k\gg k_{\crit}$ that
\begin{align*}
 a_k Q_k {\bar n}_1^k \sum_{l=k}^{\infty} \frac{1}{a_l Q_l {\bar n}_1^{l}} 
 &\simeq e^{-G(k)}\sum_{l=k}^\infty e^{G(l)}
 \simeq \int_0^\infty e^{G'(k)l}\,dl =\frac1{-G'(k)}\simeq \frac1\eps \,,
\end{align*}
and hence, from \eqref{stst}, we obtain~\eqref{nstatasymptotics}.

\subsection{Transport equation for large clusters}
\label{Ss.transport}

For clusters much larger than the critical size, exceeding the macroscopic cluster size scale $X$ that will be determined below in Section \ref{Ss.scales}, we have 
\begin{equation}\label{Jklarge}
J_k = a_k \big(n_1{-}1\big) n_k - \frac{a_kq}{k^\gamma} n_k + b_k n_k- b_{k+1}n_{k+1} \simeq a_k \eps n_k \quad\text{for } k\gg X \gg k_\crit \,.
\end{equation}
Hence, we can approximate the evolution of clusters in this regime by 
\begin{equation}\label{e:transportfeps}
\partial_{t} n(k,t)  +\eps \partial_{k}\big(k^{\alpha}n(k,t)  \big)  =-rk^\rexp n(k,t)
\end{equation}
where for large cluster sizes, we treat $k$ a continuous variable and we represent the discrete concentration $n_{k}(t)$ by a continuous
concentration $n(k,t)$.

\subsection{Monomers and the nucleation flux}
\label{Ss.fluxmatch}
The behavior of $n$ for  clusters that are much larger than the critical size, but much smaller than $X$, is given by
the quasistationary solutions depending on $n_1(t)$, that is \eqref{nstatasymptotics} implies 
\begin{equation}\label{fasym3}
k^\alpha n(k,t) \simeq \frac{J(n_1(t))}{\eps }\, \qquad \text{ for }k_\crit\ll k \ll X \,.
\end{equation}

We see from \eqref{Jasymp} that small changes in $\eps$ yield large changes in the flux $J(1+\eps)$. 
In order to obtain variations of order one during the evolution we introduce a rescaled concentration of monomers.
More precisely, we introduce for fixed $0<\eps\ll 1$ the new variable $u$ via
\[
n_{1}(t)=1+\eps+\Big(\frac{\eps}{q}\Big)^{\frac{1}{\gamma}}u(t) = \bar n_1+ \frac{u(t)}{k_\crit}\,.
\]
Note also that we have $n_{1}{-}1\simeq \eps$ at leading order as long as $u=O(1)$.
Hence, as long as $u$ remains of smaller order, we may approximate $J(n_1)$ with $n_1$ in place of $\bar n_1$ 
in the derivation of~\eqref{Jasymp} to arrive at
\begin{equation}\label{Jinftydef}
	J(n_1) \simeq J_\infty e^u \,.
\end{equation}
Indeed, the relation~\eqref{Jinftydef} follows by maximizing  
\[
\tilde G(k) = -k\log n_1 + 
\int_1^k \log(1+ql^{-\gamma})\,dl + C = G(k)-k
\log\left(\frac{n_1}{\bar n_1}\right).
\]
The maximum occurs at $\tilde k$ satisfying $n_1=1+q\tilde k^{-\gamma}$, so 
\[
\tilde k = \left(\frac q{n_1-1}\right)^{\frac 1\gamma} = 
\left(\frac q\eps\right)^{\frac1\gamma}
\left(\frac{n_1-1}\eps\right)^{-\frac1\gamma}
= k_\crit \left(1+\frac{u}{\eps k_\crit}\right)^{-\frac1\gamma} = k_\crit + O(\eps^{-1}).
\]
Hence, by following the same derivation as for~\eqref{Jasymp}, we arrive at $J(n_1)\simeq J(\bar n_1) e^{-\tilde G(\tilde k)+G(k_\crit)}$. Now $\tilde k=k_\crit (1+o(1))$, so  
\[
\tilde k
\log\left(\frac{n_1}{\bar n_1}\right) = \tilde k \log\left(1+\frac u{k_\crit \bar n_1}\right) = u + o(1)\,,
\]
while
$
G(\tilde k)-G(k_\crit) = O( G''(k_\crit)(\tilde k-k_\crit)^2) = O(k_\crit^{-1-\gamma}\eps^{-2}) = O(\eps^{\frac1\gamma-1})
$
justifying the asymptotic expansion~\eqref{Jinftydef}.

\subsection{Evolution of monomers}
\label{Ss.monomers}

We introduce a constant $L\gg1$ as a multiplicative cutoff between subcritical and supercritical clusters,
and  approximate the equation for
the monomers by
\[
\partial_{t}n_{1}=-J_{1}-\sum_{k=1}^{Lk_{\crit}}J_{k}-\sum_{k=Lk_{\crit}
+1}^{\infty}J_{k}+S \,.
\]
In the second sum we use the approximation \eqref{Jklarge} for the fluxes $J_k$ and approximate the sum by an integral
\[
\sum_{k=Lk_{\crit}+1}^{\infty}J_{k}\simeq (n_{1}{-}1)  \int
_{Lk_{\crit}+1}^{\infty}n(k,t)  k^{\alpha}\, dk \qquad\text{for } L \gg 1 . 
    \]
In \AppendixSupp~\ref{sec:quasistationary},
we justify that in the region of of subcritical clusters, the evolution follows a quasistationary law on an algebraically large time scale in $\eps^{-1}$,
which allows us to approximate $J_{k}$ by $J(n_1(t))$ for $1\leq k\leq L k_{\crit}$. 

Using also that $n_{1} - 1$ is of order $\eps$ we obtain the following equation for the concentration of
monomers
\begin{equation}\label{eq:monomerseps}
\frac1{k_\crit} {\partial_{t}u}
= -\big(  Lk_{\crit}+1\big)  J(n_1(t))-\eps\int_{Lk_{\crit}+1}^{\infty
}f(y,t)  y^{\alpha}dy+S \,,
\end{equation}
together with equation~\eqref{e:transportfeps} for the supercritical clusters and the boundary condition~\eqref{fasym3}.
We will see later that the term $(Lk_{\crit}+1)J$ in the
contribution of the monomers is negligible compared to the integral term
during the whole dynamics. We assume that for the moment and will check
that \textit{a posteriori}.

\subsection{Identification of scales}\label{Ss.scales}

We introduce new units and variables for the cluster size $X,x$ time
$T,\tau$ and cluster density $F,f$, respectively, via
\begin{equation}\label{eq:rescaling}
k=Xx\,,\qquad  t=T\tau\,,\qquad n=Ff \,.
\end{equation}
We will obtain the limit model~\eqref{eq:limit:u}--\eqref{eq:limit:boundary} by choosing the scales
\begin{equation}\label{eq:scales}
	X = \biggl( \frac{\eps}{k_\crit J_{\infty}}\biggr)^{\frac{1}{2{-}\alpha}} , \quad
	T=\frac{X^{1-\alpha}}{\eps} , \quad
	F = \frac{k_\crit}{X^2}, \quad
	S= \frac{1}{T k_\crit}, \quad
	R = \frac{\rfac}{T X\rexp},
\end{equation}
with $k_\crit$ as defined in~\eqref{criticalsize} and $J_\infty$ given in~\eqref{Jasymp}. 
We have also included the size of the monomer source $S$ in this formula and introduce a rescaled depletion rate $\rfac >0$.
Notice that from \eqref{Jasymp}, $J_\infty$ is exponentially small in the parameter $\eps$, which implies similar exponential dependencies for all scales $X, T,F$, and $S$. 

Using the scales in \eqref{eq:rescaling},
the monomer equation~\eqref{eq:monomerseps} becomes 
\begin{equation}\label{eq:limit:scaled:u}
\partial_{\tau}u=-\big(  Lk_{\crit}+1\big)  J(n_1(\tau))T k_\crit-\eps X^{1+\alpha}FT k_\crit\int_{\frac{Lk_{\crit}
+1}{X}}^{\infty}f(x,\tau)  x^{\alpha}dx+ST k_\crit \,,
\end{equation}
and similarly the transport equation for the large clusters~\eqref{e:transportfeps} together with its boundary conditions~\eqref{fasym3} take the form
\begin{align}
\partial_{\tau}f(x,\tau)  +\frac{T\eps}{X^{1-\alpha}}\partial_{x}\big(  x^{\alpha}f(x,\tau)  \big)   &  =-\big(TX^\rexp R\big)  x^\rexp f(x,\tau) \label{eq:limit:scaled:transport}\\
x^\alpha f(x,\tau)   &\simeq\frac{J(n_1(\tau))}{F\eps X^{\alpha}} ,
\qquad \text{for } \frac{k_\crit}{X} \ll x \ll 1\,. \label{eq:limit:scaled:boundary}
\end{align}

We recall from \eqref{Jinftydef} that $J(n_1(\tau))=J_{\infty}e^{u(\tau)}$, where  $J_{\infty}$ is an exponentially small quantity defined in~\eqref{Jasymp}
that gives the order of magnitude of the fluxes through the critical size.
As explained in Section~\ref{Ss.description}, the nucleation rate, given by the flux through the critical size manifested as the boundary condition in~\eqref{eq:limit:scaled:boundary}, and the condensation rate, given by the coagulation rate of monomers with macroscopic clusters manifested as the scale of the integral in~\eqref{eq:limit:scaled:u}, need to be of the same order of magnitude. It is readily checked that~\eqref{eq:scales} implies
\[
\frac{J_{\infty}}{F\eps X^{\alpha}}    =1
\qquad \text{ and } \qquad \eps X^{1+\alpha}FT k_\crit=1\,.
\]
The time scale that yields growth of the
macroscopic clusters due to the transport on the left hand side of~\eqref{eq:limit:scaled:transport} in that scale is given by ${T\eps}/{X^{1-\alpha}}$, which becomes equal to one by the choice~\eqref{eq:scales}.

Notice that the choice of $S$ in~\eqref{eq:scales} implies that the changes of $u$ in~\eqref{eq:limit:scaled:u} are of order one. This justifies the expansion of the flux in~\eqref{Jinftydef} \textit{a posteriori}
This formula yields the relation between $\eps$ and the source $S$.
A different choice of $S$ might still result in oscillatory behavior, with an amplitude change in $u$ not of order $O(1)$, although the expansion of the flux in~\eqref{Jinftydef} is valid for $u$ up to order~$o(\eps^{-\frac{1}{\gamma}})$.

In order to conclude the derivation of the limit model~\eqref{eq:limit:u}--\eqref{eq:limit:boundary}, it only remains to justify neglecting the fluxes in the subcritical region, or equivalently to show that the
term $(  Lk_{\crit}+1)  J(n_1(\tau))T k_\crit$ in~\eqref{eq:limit:scaled:u} is negligible.
The contribution $  Lk_{\crit}+1  $ is algebraically large in $\eps^{-1}$.
We estimate the exponential terms contained in the product $J(n_1(\tau))T$,
which are 
\begin{equation}\label{e:negligible_fluxes}
J_{\infty}\Big( \frac{1}{J_{\infty}}\Big)^{\frac{1-\alpha}{2-\alpha}}=J_{\infty}^{1-\frac{1-\alpha}{2-\alpha}}=J_{\infty}^{\frac{1}{2-\alpha}}\,.
\end{equation}
Since $\frac{1}{2-\alpha}>0$, it follows that this term is exponentially
small. Therefore the flux terms due to subcritical particles yield a
negligible contribution and~\eqref{eq:limit:scaled:u} can be formally approximated by~\eqref{eq:limit:u}.

\section{Periodic solutions for the limit problem}\label{S.hopf}

Our goal is to argue, by formal means, that steady state solutions of the limit model~\eqref{eq:limit:u}--\eqref{eq:limit:boundary}
undergo a Hopf bifurcation as the rescaled removal parameter $\rfac$ is varied, 
for a wide range of values of $\alpha$ and $b$.

To our knowledge, a rigorous Hopf bifurcation theorem has not yet been proven for a model of this type,
though it seems plausible that one might extend existing methods for retarded functional differential
equations with finite delays (RFDE) which are based upon 
rescaling to fix the temporal period and Lyapunov-Schmidt reduction
\cite{Hale_book,Diekmann_book}.  

It is well known, that in order to prove the existence of periodic solution by means of a Hopf bifurcation, one has to show:
\begin{enumerate}
	\item existence of steady states for a family of dynamical systems parametrized by $\theta\in \Theta$;
	\item the linear stability of these steady states changes at a critical value $\theta_\crit$ due to fact that two complex conjugated eigenvalues of the linearized system cross from $\{\re(\lambda)<0\}$ to $\{\re(\lambda) > 0 \}$;
	\item for generic dynamic systems, the periodic solutions exists for small $|\theta - \theta_\crit|$ either in the parameter region $\{ \theta < \theta_\crit  \}$ or $\{\theta > \theta_\crit \}$.
	These periodic solutions are stable if the steady states found in $1.$ are unstable for the range of $\theta$, where the periodic solutions exist.
\end{enumerate}

We carry out the three steps of this program for the limit model~\eqref{eq:limit:u}--\eqref{eq:limit:boundary} in Subsections~\ref{Ss.steadyfirst}, \ref{Ss.crossing}, \ref{Ss.supercritical}, respectively.
To simplify the Hopf bifurcation analysis, it is convenient to introduce a new size variable $z$ and relabeled flux ${h}(z,\tau)$ defined by
\begin{equation}\label{zparms}
	z = x^{1-\alpha} \qquad\text{and}\qquad 
	h(z,\tau) = \frac{x^\alpha g(x,\tau)}{1-\alpha} \,,
\end{equation}
In those variables the system~\eqref{eq:limit:u}--\eqref{eq:limit:boundary} takes the form
\begin{align}
	\partial_{\tau}u  &  =1-\int_{0}^{\infty}h(z,\tau)  z^{\nu}\,dz\,,\label{ueq3}\\
	\partial_{\tau}h(z,\tau)  +\partial_{z}  h(z,\tau )   &  =-\rfac z^{\beta}h(z,\tau)
	\,, \qquad z>0 \,, \label{transporth}\\
	h(  0,\tau)   &  =e^{u(\tau)}\,,\label{hinitial}
\end{align}
where the exponents $\beta$ and $\nu$ are nonnegative and given by 
\begin{equation}\label{d:nubeta} 
\beta = \frac{\rexp}{1-\alpha}\,,
\qquad
\nu = \frac{\alpha}{1-\alpha}\,. 
\end{equation}

For any solution of this system, such as a steady-state or time-periodic solution,
which exists for all times $\tau\in\R$, one finds,
by integrating \eqref{transporth} along characteristics emerging from $z=0$,
that necessarily
\begin{equation}\label{hform}
h(z,\tau)  =\exp(u(\tau{-}z))
\exp\Bigl(  -\frac{\rfac}{\beta+1}z^{\beta+1}\Bigr) \,.
\end{equation}
Using this expression in  \eqref{ueq3} we find that $u$ must satisfy an RFDE with infinite delay, namely
\begin{equation}\label{e:u:closed}
\partial_{\tau}u=1-\int_{0}^{\infty}e^{u(  \tau{-}z)  }\exp\Bigl(
-\frac{\rfac}{\beta+1} z^{\beta+1}\Bigr)  z^{\nu}dz \,.
\end{equation}

\subsection{Steady states and their stability}\label{Ss.steadyfirst}
For any $\rfac >0$ we have a constant solution $u\equiv u_0$ of this equation, given by
\begin{equation}\label{e:u0}
1 =e^{u_0}
\int_{0}^{\infty}\exp\Bigl(  -\frac{\rfac}{\beta+1}z^{\beta+1} \Bigr)  z^{\nu}dz \, .
\end{equation}
The corresponding steady state $h=h_0(z)$ is then given by \eqref{hform} with $u$ replaced by $u_0$.

Next we consider the linear stability of these steady states as solutions of \eqref{ueq3}--\eqref{hinitial}.
After linearizing about $u_0$, $h_0$, 
we find nonzero solutions proportional to $e^{\lambda\tau}$ exist, $\lambda\in\C$, if and only if 
\begin{equation}\label{e:eval1}
\lambda = - e^{u_0}
\int_{0}^{\infty}\exp\Bigl( -\frac{\rfac}{\beta+1}z^{\beta+1} -\lambda z\Bigr)  z^{\nu}dz\,.
\end{equation}
It is convenient to study this equation after the rescaling
$\rfac^{-\frac1{\beta+1}}\lambda \mapsto \lambda$. With
\begin{equation}\label{d:G}
G_{\beta,\nu}(\lambda) := \int_{0}^{\infty}\exp\Bigl(  -\frac 1{\beta+1} {z^{\beta+1}}-\lambda z\Bigr)  z^{\nu}dz
\,,
\end{equation}
the eigenvalue equation \eqref{e:eval1} then takes the form
\begin{equation}\label{e:eval2}
\vth\lambda + G_{\beta,\nu}(\lambda) = 0,
\qquad\text{where}\quad
\vth = \rfac^{\frac1{\beta+1}}  G_{\beta,\nu}(0)\,.
\end{equation}
Note $\lambda=0$ is never an eigenvalue since $G_{\beta,\nu}(0)>0$.
Also note $G_{\beta,\nu}'(\lambda)=-G_{\beta,\nu+1}(\lambda)$, 
\begin{equation}\label{Gest1}
G_{\beta,\nu}(0) \geq |G_{\beta,\nu}(\lambda)| 
\quad\text{and}\quad
|G_{\beta,\nu}'(0)| \geq |G_{\beta,\nu}'(\lambda)| 
\quad\text{whenever  $\re\lambda\ge0$.}
\end{equation}
It follows that we have linear stability for large values of the parameter $\vth$: 
\begin{lemma} 
There are no eigenvalues in the closed right half plane
if $\vth>|G_{\beta,\nu}'(0)|$, i.e.,
\[
\eta > \left(\frac{G_{\beta,\nu+1}(0)}{G_{\beta,\nu}(0)}\right)^{\beta+1}.
\]
\end{lemma}
The reason is that if $\re\lambda\ge0$, then with $G=G_{\beta,\nu}$ we have
\[
\vth > |G'(0)| \geq \re\left( \frac{G(0)-G(\lambda)}\lambda\right) \geq -\re \frac{G(\lambda)}\lambda.
\]

\subsection{Eigenvalue crossings}\label{Ss.crossing}
A Hopf bifurcation should occur provided that some branch of solutions $\lambda=\lambda(\vth)$
of \eqref{e:eval2} crosses the imaginary axis transversely. 
This means that $\re\lambda=0$ and $\re \frac{d\lambda}{d\vth}\ne0$ 
at some particular $\vth=\vth_0$.  
Since $(\vth+G'(\lambda))\frac{d\lambda}{d\vth}=-\lambda$ along the branch,
we find
\begin{equation}\label{e:dlam}
\sign\re \frac{d\lambda}{d\vth} = \sign\re\Bigl(-\lambda\bigl(\vth+\overline{G'(\lambda)}\bigr)\Bigr) 
= \sign \frac{d}{dt} \re G(it) = \sign 
\frac{d}{dt}\arg G(it),
\end{equation}
if $\lambda=it$ with $t>0$.  
Thus the criteria for Hopf bifurcation become the following: First,
there should exist $t_0>0$ such that for $t=t_0$,
\begin{equation}\label{c:hopf1}
\re G_{\beta,\nu}(it)=0 
\quad\text{and}\quad
\im G_{\beta,\nu}(it)<0 ,
\quad\text{or equivalently}\quad \arg G_{\beta,\nu}(it) = -\frac\pi2 \ (\mbox{mod}\,2\pi)\,.
\end{equation}
This is necessary and sufficient for \eqref{e:eval2} to hold with $\vth=-G(it)/it>0$. 
Second, the transversality condition holds if and only if 
\begin{equation}\label{c:hopf2}
\frac{d}{dt} \re G_{\beta,\nu}(it) \ne0, 
\qquad\text{or equivalently }\quad \frac{d}{dt}\arg G_{\beta,\nu}(it) \ne0.
\end{equation}
Thus we can provide evidence that a Hopf bifurcation occurs (for some $\vth$)
and infer the direction of eigenvalue crossings by identifying zero crossings on 
the graph of $\frac\pi2+\arg G(it)$ (mod $2\pi$).
In the original time scale $\tau$ of the model,
such zeros correspond to oscillations with wave number
\begin{equation}\label{d:waveno}
\kappa = t\rfac^{\frac1{\beta+1}} = \frac{t\vth}{G_{\beta,\nu}(0)}\,, 
\qquad G_{\beta,\nu}(0) = (\beta+1)^{\frac{\nu-\beta}{\beta+1}} \Gamma\left(\frac{\nu+1}{\beta+1}\right).
\end{equation}

\subsection{Direction of bifurcation}\label{Ss.supercritical}

In this section we identify computable criteria that should determine
the direction of bifurcation (i.e., whether small time-periodic solutions appear for
$\rfac>\rfac_0$ or $\rfac<\rfac_0$). 
We posit bifurcating solutions have variable period $2\pi/\kappa$
and scale time using the variable wave number $\kappa$ near a value $\kappa_0$ given by \eqref{d:waveno} 
at a putative bifurcation point $t=t_0$.
Rewriting~\eqref{e:u:closed} in terms of the 
constant solution $u_0$ from \eqref{e:u0} and rescaled variables given by
\begin{equation}
u(\tau) = u_0 + U(s) ,\qquad s = \kappa\tau, \qquad y=\kappa z,
\end{equation}
we find that the $2\pi$-periodic function $U$, constant $\kappa$ and bifurcation parameter $\rfac$ should satisfy 
\begin{equation}
\kappa \partial_s U  + \int_0^\infty \left(e^{U(s-y)} -1\right)
\exp\left(\frac{-\rfac}{\beta+1}\frac{y^{\beta+1}}{\kappa^{\beta+1}} \right)y^\nu\,dy 
\frac{e^{u_0}}{\kappa^{\nu+1}} = 0.
\end{equation}

It will be convenient to work with the variables and parameters given by 
\begin{equation}\label{e:def:delta:mu}
v(s) = e^{U(s)}-1,\qquad \delta = \frac{\kappa^{\nu+2}}{e^{u_0}}, \qquad \mu = \frac{\rfac}{(\beta+1)\kappa^{\beta+1}},
\qquad
g(\mu,y)  = y^\nu \exp(-\mu y^{\beta+1}).
\end{equation}
Thus we seek a $2\pi$-periodic function $v$, constant $\delta$ and bifurcation parameter $\mu$ so 
\begin{equation}\label{v:bif}
\delta\partial_s v + (1+v) (\gamma * v) = 0,
\qquad\text{where}\quad
(\gamma * v)(s) = \int_0^\infty g(\mu,y)v(s-y)\,dy. 
\end{equation}
The function $v\equiv0$ is a trivial solution for any $\delta$ and $\mu$, 
and the problem is equivariant with respect to translation in $s$.
We expect bifurcation from this branch for $\delta$, $\mu$
respectively near values $\delta_0$, $\mu_0$ coming from the 
parameters $t_0$, $\rfac_0$ by the formulas above.

In terms of an amplitude parameter $\eps$, we seek formal expansions
\[
v = \eps v_1+\eps^2 v_2+\ldots,
\quad 
\delta = \delta_0+\eps \delta_1 + \eps^2 \delta_2 + \ldots,
\quad
\mu = \mu_0+\eps \mu_1 + \eps^2 \mu_2 + \ldots.
\]
We will find $\delta_1=\mu_1=0$ as is typical for Hopf bifurcation.
Provided certain nondegeneracy conditions hold, 
we find the quantities $\delta_2$ and $\mu_2$ can be expressed in terms of the quantities
\begin{equation}\label{e:def:hatg}
\hat g_j(k) = \int_0^\infty e^{-iky} g_j(y)\,dy,
\qquad 
g_j(y) = \partial_\mu^j g(\mu_0,y) = (-y^{\beta+1})^jy^\nu\exp(-\mu_0 y^{\beta+1}). 
\end{equation}
In fact, our computations will show that $\mu_2$ and $\delta_2$ are determined by the relations 
\begin{equation}\label{mu2delta2}
  \mu_2 \re \hat g_1(1) =  
  \frac{-\delta_0^2\re \hat g_0(2)}{4|\hat g_0(2)+2i\delta_0|^2} , \qquad
\delta_2 + \mu_2\im \hat g_1(1) = 
 -\frac{\delta_0}4 + \frac{\delta_0^2(2\delta_0 + \im \hat g_0(2))}
{4|\hat g_0(2)+2i\delta_0|^2},
\end{equation}
provided the transversal crossing condition and the nonresonance condition
$\hat g_0(2)+2i\delta_0\ne0$ hold.

Returning to the original parameters 
$\rfac=\rfac_0(1+\eps^2\rfac_2+\ldots)$, 
$\kappa=\kappa_0(1+\eps^2\kappa_2+\ldots)$, and taking from~\eqref{e:u0} into account that $e^{u_0} G_{\beta,\nu}(0)= \rfac^{\frac{\nu+1}{\beta+1}}$,
we find 
\begin{equation}\label{r:theta2}
\rfac_2
=  (\nu+2)\frac{\mu_2}{\mu_0} + (\beta+1)\frac{\delta_2}{\delta_0}
,\qquad
\kappa_2
= \frac{\nu+1}{\beta+1}\frac{\mu_2}{\mu_0} + \frac{\delta_2}{\delta_0}.
\end{equation}
The sign of $\rfac_2$ determines the direction of the bifurcation. 
It is plausible that this determines the stability of the bifurcating periodic solutions in a similar way as for ODEs, 
and RFDEs with finite delays~\cite{Diekmann_book}.  
Namely, if the constant flux solution with $u\equiv u_0$ is stable for $\rfac$ on one side of the bifurcation point $\rfac_0$,
then bifurcating periodic solutions are stable if they appear for $\rfac$ on the opposite side, and unstable otherwise.
As the constant-flux solution is stable for $\rfac$ sufficient large, we therefore expect that {\em at the largest bifurcation point $\rfac_0$,
the new branch of periodic solutions is stable provided that $\rfac_2$ is negative.} 
We call this case \emph{supercritical} and will provide numerical evidence that indeed such bifurcations occur in Section~\ref{Ss.numerical}.

In the remainder of this section we derive the relations in \eqref{mu2delta2}.
By an integration by parts starting from \eqref{e:dlam}, we have $\lambda=it=-\bar\lambda$, 
since $\frac{d}{dy}e^{-ity}=-it e^{-ity}$. Moreover, $G'_{\beta,\nu}(\lambda)=-G_{\beta,\nu+1}(\lambda)$ and hence by using \eqref{e:eval2}, we obtain
\begin{align*}
-\bar\lambda G'_{\beta,\nu}(\lambda) &= 
 -it G_{\beta,\nu+1}(it)
= \int_0^\infty\left(\frac{d}{dy} e^{-ity} \right)
 y^{\nu+1} \exp\Bigl(  -\frac{y^{\beta+1}}{\beta+1}\Bigr)  dy
\\ & = - \int_0^\infty e^{-ity}( (\nu+1)y^\nu - y^{\nu+1+\beta})
 \exp\Bigl(  -\frac{y^{\beta+1}}{\beta+1}\Bigr)  dy
\\ &  = -(\nu+1)G_{\beta,\nu}(it) + G_{\beta,\nu+1+\beta}(it)
\\& = (\nu+1)(it)\vth - \frac {\hat g_1(1)} {t^{\nu+2+\beta}}  . 
\end{align*}
Here, we used~\eqref{e:def:hatg} and thus arrive at
\begin{equation}\label{e:deriv:crossing}
\sign\re\frac{d\lambda}{d\vth} 
 = - \sign\re\hat g_1(1) \,.
\end{equation}
We require this quantity is nonzero, which is equivalent to transversal eigenvalue crossing. 

We proceed to compute. Expanding $\gamma = \gamma_0+ \eps \gamma_1+\eps^2\gamma_2 +\ldots$ we find
$\gamma_0 = g_0$, $\gamma_1 = \mu_1 g_1$, $\gamma_2 = \mu_2 g_1+\frac12\mu_1^2 g_2$.
Plugging into \eqref{v:bif}, at the respective orders $\eps$, $\eps^2$, and  $\eps^3$ we find
\begin{align}
0 &= L v_1 := \delta_0\partial_s v_1 + \gamma_0*v_1 ,
\label{b:eps1}\\
0 &= L v_2 + 
(\delta_1\partial_s+\gamma_1*)v_1 
+ v_1 \gamma_0*v_1 ,
\label{b:eps2}\\
0 &= L v_3 + 
(\delta_2\partial_s+\gamma_2*)v_1 
+ v_2(\gamma_0*v_1)+v_1(\gamma_0*v_2)
\label{b:eps3} \\
& \qquad + (\delta_1\partial_s+\gamma_1*)v_2 + v_1(\gamma_1*v_1) .
\nonumber 
\end{align}
By what we have said in subsection~\ref{Ss.crossing}, equation \eqref{b:eps1}
should have two-dimensional kernel $N_0$ spanned by $\cos s$ and $\sin s$ (or $e^{\pm is}$). 
Due to translation equivariance of the problem, we may suppose that the amplitude and phase 
of this mode of $v$ are normalized so that $v_1(s) = \cos s$.

Now, in terms of Fourier series $v = \sum_{k\in\Z}\hat v(k)e^{iks}$, an equation $Lv=w$ 
corresponds to 
\[
\hat L(k)\hat v(k)=\hat w(k) \quad\text{for all $k\in\Z$,}
\qquad \hat L(k) = ik\delta_0+\hat g_0(k). 
\]
This is solvable if and only if the {nonresonance condition} holds, 
namely $\hat L(k)\ne0$ for all $k$ such that $\hat w(k)\ne0$.
Since $\hat L(\pm1)=0$ by assuming the first criteria~\eqref{c:hopf1} for a Hopf transition, we always require $\hat w(\pm1)=0$.

Moreover, for the following formal analysis we require the nonresonance condition
\begin{equation}\label{e:nonresonance}
 \hat L(k) = \overline{\hat L(-k)} = ik\delta_0+\hat g_0(k)\ne 0\qquad\text{holds for } k=\pm2 . 
\end{equation}
Applying these considerations to \eqref{b:eps2}, since $L(\cos s)=0$ we see that the term 
\[
v_1(\gamma_0*v_1)=(\cos s)(\delta_0 \sin s) = -\frac{i\delta_0}4 e^{2is}+c.c.,
\]
is nonresonant, provided $\hat L(\pm2)\ne0$ as we have assumed.
However the terms in $\delta_1$ and $\mu_1$ are resonant. Since $\cos s = \frac12 e^{is}+c.c.$ we find
these terms yield
$
(i\delta_1+ \mu_1\hat g_1(1))\frac12 e^{is} + c.c.
$
Since we presume $\re\hat g_1(1)\ne0$ from~\eqref{c:hopf2},
solvability of \eqref{b:eps2} requires $\mu_1=0$, and also $\delta_1=0$. The solution for $v_2$ is 
\begin{equation}
v_2(s) = \hat v_2(2) e^{2is}+c.c., \qquad \hat v_2(2) = \frac{i\delta_0}4\frac{1}{2i\delta_0+\hat g_0(2)}.
\end{equation}

Next, we consider the solvability of \eqref{b:eps3}.  The terms involving $\delta_1$ and $\gamma_1$ vanish.
It remains to show that the resonant terms are removed by a particular choice of $\delta_2$ and $\mu_2$.
For this purpose, noting $\hat g_0(1)=-i\delta_0$ we observe
\[
\gamma_0*v_1 = \hat g_0(1)\hat v_1(1)e^{is}+c.c. = -\frac{i\delta_0}2 e^{is}+c.c., \qquad 
\gamma_0*v_2 = \hat g_0(2)\hat v_2(2) e^{2is}+c.c.,
\]
therefore the resonant terms on the right hand side of \eqref{b:eps3} are 
\begin{align*}
& e^{is}\Bigl(
(\delta_2 i+ \mu_2\hat g_1(1))\hat v_1(1)
+ \hat v_2(2)\overline{\hat g_0(1)\hat v_1(1)}
+ \overline{\hat v_1(1)}\hat g_0(2)\hat v_2(2)
\Bigr) + c.c.
\\
& =\ \ 
\frac12 e^{is}\left( \delta_2 i +\mu_2\hat g_1(1) + (i\delta_0+\hat g_0(2)) \hat v_2(2) \right) + c.c.,
\end{align*}
Since 
$i\delta_0+\hat g_0(2) = \hat L(2)-i\delta_0$,
the resonant terms vanish if $\mu_2$ and $\delta_2$ are given by 
\begin{align}\label{mu2form}
  \mu_2 \re \hat g_1(1) &=  
- \re \left( \frac{i\delta_0}4 \frac{\hat L(2)-i\delta_0}{\hat L(2)}\right)
 = \frac{-\delta_0^2}{4|\hat L(2)|^2} \re \hat g_0(2),
 \\
\label{delta2form}
\delta_2 + \mu_2\im \hat g_1(1) &= 
- \im \left( \frac{i\delta_0}4 \frac{\hat L(2)-i\delta_0}{\hat L(2)}\right)
= -\frac{\delta_0}4 + \frac{\delta_0^2}{4|\hat L(2)|^2} (2\delta_0 + \im \hat g_0(2)).
\end{align}

\section{Ranges of parameters yielding bifurcation}\label{s:parms}
In this section we provide some numerical and analytical information regarding the 
occurrence and direction of Hopf bifurcations as one varies the removal factor $\rfac$
in the limit model~\eqref{ueq3}--\eqref{hinitial},
for a variety of cases involving the exponents $\beta=r/(1-\alpha)$ and $\nu=\alpha/(1-\alpha)$ from \eqref{d:nubeta}.  
A particularly intriguing finding is that an infinite number of bifurcation points appear
when $\nu=\alpha=0$ and $\beta=r$ is an odd integer greater than $1$.

\subsection{The case \texorpdfstring{$\beta=r=0$}{β=0}, \texorpdfstring{$\nu\ge0$}{ν≥0}}\label{Ss.beta0}
We have explicit formulas in this case, since whenever $\re\lambda>-1$, 
the substitution $(1+\lambda)y=s$ and a contour deformation argument shows
\[
G_{0,\nu}(\lambda)=\frac{\Gamma(\nu+1)}{(1+\lambda)^{\nu+1}}.
\]
This case corresponds to removal of clusters in \eqref{nkeq} at rate independent of size, since $r=0$. 
When $\nu$ is an integer, equation~\eqref{e:u:closed} can be reduced to an ODE of some order 
by repeated differentiation. 
We comment on the connection to similar models in the literature in the Appendix~\ref{Ss.momentmodels}.

In the present case, we can rigorously describe all eigenvalue crossings and prove transversality.
\begin{proposition}
In case $\beta=0$ and $\nu\ge0$, an eigenvalue $\lambda=it$ 
with $t>0$ occurs for some $\rfac>0$ if and only if 
\[
\tan^{-1} t = \omega_k := \frac{\pi}2\frac{1+4k}{1+\nu} 
\quad\mbox{for some integer $k$ satisfying $0\le 4k<\nu$.}
\]
The transversal crossing condition~\eqref{c:hopf2} holds for all $k$, 
with $\sign \re \frac{d\lambda}{d\vth} <0$.
\end{proposition}
The reason this is true is that the condition \eqref{c:hopf1} for 
eigenvalue crossings takes the form 
\[
-{\arg G_{0,\nu}(it)} = (\nu+1)\arg (1+it) = (\nu+1)\tan^{-1} t = \frac\pi2\ (\mbox{mod}\,2\pi), 
\]
and transversality follows from \eqref{e:dlam}.
We note that large $\rfac$ (or $\vth$) corresponds to small $t$, 
so that as $\rfac$ decreases from large values (where one has linear stability), a first crossing 
appears at
\begin{equation}\label{e:beta0val}
t_0 = \tan \omega_0, \qquad \rfac_0 = \frac{\vth_0}{\Gamma(\nu+1)} = \frac{\cos^{\nu+2}\omega_0}{\sin\omega_0},
\qquad \kappa_0 = \cos^{\nu+1}\omega_0.
\end{equation}

In summary, for $\beta=0$ we find that whenever $\nu>0$, bifurcation from a stable state
should occur at wave number $\kappa_0$ as $\rfac$ decreases through the value in \eqref{e:beta0val}.
However, no bifurcation occurs for $\nu=0$.
The values of $\rfac_2$ and $\kappa_2$ are in principle explicitly computable, but the expressions become quickly impractical.

\subsection{Numerical bifurcation curves}\label{Ss.numerical}

\begin{figure}%[ht]
\centering
\begin{tikzpicture}
\begin{axis}[colormap/viridis,cycle list={[colors of colormap={200,400,600,800,1000}]},
width=0.47\textwidth, axis lines=left, grid=major, ymin=0, xmin=0,xmax=10, xlabel={$\ \beta$}, restrict y to domain=0:1, ylabel={$\rfac_0$},
every axis x label/.style={at={(current axis.right of origin)},anchor=south,},
every axis y label/.style={at={(current axis.north west)},anchor=north east,}]
\addplot+[thick,mark=] table [x=beta,y=theta0_1]{BDoscillate-theta-kappa.csv};\label{p:theta0_1}
% \addplot+[thick,mark=] table [x=beta,y=theta0_2]{BDoscillate-theta-kappa.csv};\label{p:theta0_2}
\addplot+[thick,mark=] table [x=beta,y=theta0_3]{BDoscillate-theta-kappa.csv};\label{p:theta0_3}
% \addplot+[thick,mark=] table [x=beta,y=theta0_4]{BDoscillate-theta-kappa.csv};\label{p:theta0_4}
\addplot+[thick,mark=] table [x=beta,y=theta0_5]{BDoscillate-theta-kappa.csv};\label{p:theta0_5}
% \addplot+[thick,mark=] table [x=beta,y=theta0_6]{BDoscillate-theta-kappa.csv};\label{p:theta0_6}
\addplot+[thick,mark=] table [x=beta,y=theta0_7]{BDoscillate-theta-kappa.csv};\label{p:theta0_7}
% \addplot+[thick,mark=] table [x=beta,y=theta0_8]{BDoscillate-theta-kappa.csv};\label{p:theta0_8}
\addplot+[thick,mark=] table [x=beta,y=theta0_9]{BDoscillate-theta-kappa.csv};\label{p:theta0_9}
\end{axis}
\end{tikzpicture}\hfill%
\begin{tikzpicture}
\begin{axis}[colormap/viridis,cycle list={[colors of colormap={200,400,600,800,1000}]},
width=0.47\textwidth, axis lines = left, grid=major, ymin=0, xmin=0,xmax=10,xlabel={$\ \beta$}, ylabel={$\kappa_0$}, 
every axis x label/.style={at={(current axis.right of origin)},anchor=south,},
every axis y label/.style={at={(current axis.north west)},anchor=north east,}]
\addplot+[thick,mark=] table [x=beta,y=kappa0_1]{BDoscillate-theta-kappa.csv};\label{p:kappa0_1}
%\addplot+[thick,mark=] table [x=beta,y=kappa0_2]{BDoscillate-theta-kappa.csv};\label{p:kappa0_2}
\addplot+[thick,mark=] table [x=beta,y=kappa0_3]{BDoscillate-theta-kappa.csv};\label{p:kappa0_3}
%\addplot+[thick,mark=] table [x=beta,y=kappa0_4]{BDoscillate-theta-kappa.csv};\label{p:kappa0_4}
\addplot+[thick,mark=] table [x=beta,y=kappa0_5]{BDoscillate-theta-kappa.csv};\label{p:kappa0_5}
%\addplot+[thick,mark=] table [x=beta,y=kappa0_6]{BDoscillate-theta-kappa.csv};%\label{p:kappa0_6}
\addplot+[thick,mark=] table [x=beta,y=kappa0_7]{BDoscillate-theta-kappa.csv};\label{p:kappa0_7}
%\addplot+[thick,mark=] table [x=beta,y=kappa0_8]{BDoscillate-theta-kappa.csv};\label{p:kappa0_8}
\addplot+[thick,mark=] table [x=beta,y=kappa0_9]{BDoscillate-theta-kappa.csv};\label{p:kappa0_9}
\end{axis}
\end{tikzpicture}
\caption{$\rfac_0$ (left) and $\kappa_0$ (right) vs $\beta$ for $\nu=\frac19$~(\ref{p:theta0_1}), $\frac37$~(\ref{p:theta0_3}), $1$~(\ref{p:theta0_5}), $\frac73$~(\ref{p:theta0_7}), and $9$~(\ref{p:theta0_9}).}\label{Fig:theta0-kappa0}
\end{figure}
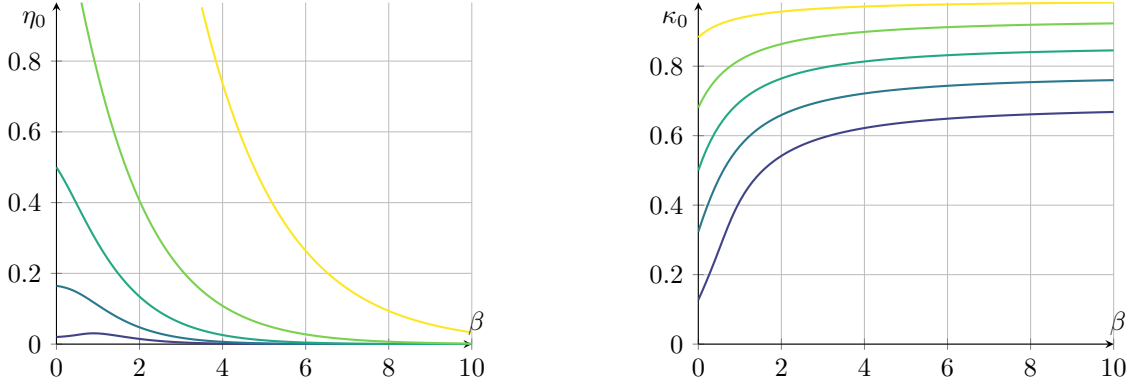

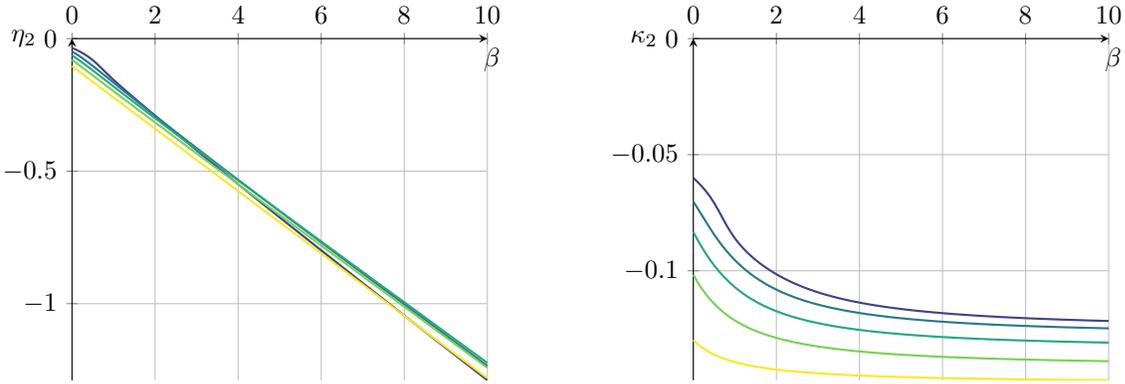
\begin{figure}%[ht]
\centering
\begin{tikzpicture}
\begin{axis}[colormap/viridis,cycle list={[colors of colormap={200,400,600,800,1000}]},
width=0.47\textwidth, axis y line=left, axis x line=top, grid=major, ymax=0, xmin=0,xmax=10,xlabel={$\ \beta$}, ylabel={$\rfac_2\ \ \ $},
every axis x label/.style={at={(current axis.right of origin)},anchor=north,},
every axis y label/.style={at={(current axis.left of origin)},anchor=east,}]
\addplot+[thick,mark=] table [x=beta,y=theta2_1]{BDoscillate-theta-kappa.csv};\label{p:theta2_1}
%\addplot+[thick,mark=] table [x=beta,y=theta2_2]{BDoscillate-theta-kappa.csv};\label{p:theta2_2}
\addplot+[thick,mark=] table [x=beta,y=theta2_3]{BDoscillate-theta-kappa.csv};\label{p:theta2_3}
%\addplot+[thick,mark=] table [x=beta,y=theta2_4]{BDoscillate-theta-kappa.csv};\label{p:theta2_4}
\addplot+[thick,mark=] table [x=beta,y=theta2_5]{BDoscillate-theta-kappa.csv};\label{p:theta2_5}
%\addplot+[thick,mark=] table [x=beta,y=theta2_6]{BDoscillate-theta-kappa.csv};%\label{p:theta2_6}
\addplot+[thick,mark=] table [x=beta,y=theta2_7]{BDoscillate-theta-kappa.csv};\label{p:theta2_7}
%\addplot+[thick,mark=] table [x=beta,y=theta2_8]{BDoscillate-theta-kappa.csv};\label{p:theta2_8}
\addplot+[thick,mark=] table [x=beta,y=theta2_9]{BDoscillate-theta-kappa.csv};\label{p:theta2_9}
\end{axis}
\end{tikzpicture}\hfill%
\begin{tikzpicture}
\begin{axis}[colormap/viridis,cycle list={[colors of colormap={200,400,600,800,1000}]},
width=0.47\textwidth, axis y line=left, axis x line=top, grid=major, ymax=0, xmin=0,xmax=10,xlabel={$\ \beta$}, ylabel={$\kappa_2\ \ \ $},
y tick label style={ /pgf/number format/.cd, fixed, precision=2, /tikz/.cd}, scaled ticks=false,
every axis x label/.style={at={(current axis.right of origin)},anchor=north,},
every axis y label/.style={at={(current axis.left of origin)},anchor=east,}]]
\addplot+[thick,mark=] table [x=beta,y=kappa2_1]{BDoscillate-theta-kappa.csv};\label{p:kappa2_1}
%\addplot+[thick,mark=] table [x=beta,y=kappa2_2]{BDoscillate-theta-kappa.csv};\label{p:kappa2_2}
\addplot+[thick,mark=] table [x=beta,y=kappa2_3]{BDoscillate-theta-kappa.csv};\label{p:kappa2_3}
%\addplot+[thick,mark=] table [x=beta,y=kappa2_4]{BDoscillate-theta-kappa.csv};\label{p:kappa2_4}
\addplot+[thick,mark=] table [x=beta,y=kappa2_5]{BDoscillate-theta-kappa.csv};\label{p:kappa2_5}
%\addplot+[thick,mark=] table [x=beta,y=kappa2_6]{BDoscillate-theta-kappa.csv};%\label{p:kappa2_6}
\addplot+[thick,mark=] table [x=beta,y=kappa2_7]{BDoscillate-theta-kappa.csv};\label{p:kappa2_7}
%\addplot+[thick,mark=] table [x=beta,y=kappa2_8]{BDoscillate-theta-kappa.csv};\label{p:kappa2_8}
\addplot+[thick,mark=] table [x=beta,y=kappa2_9]{BDoscillate-theta-kappa.csv};\label{p:kappa2_9}
\end{axis}
\end{tikzpicture}
\caption{$\rfac_2$ (left) and $\kappa_2$ (right) vs $\beta$ for $\nu=\frac19$~(\ref{p:theta2_1}), $\frac37$~(\ref{p:theta2_3}), $1$~(\ref{p:theta2_5}), $\frac73$~(\ref{p:theta2_7}), and $9$~(\ref{p:theta2_9}).}\label{Fig:theta2-kappa2}
\end{figure}

In Figure~\ref{Fig:theta0-kappa0}, we plot the largest value of $\rfac_0$ corresponding to a solution of
the bifurcation criteria~\eqref{c:hopf1} against the parameter $\beta$ for various values of $\nu$. 
The transversality condition~\eqref{c:hopf2} and nonresonance condition~\eqref{e:nonresonance} are verified numerically in 
\AppendixSupp~\ref{supp:num-ver} 
(see Table~\ref{table:transversality} and Table~\ref{table:nonresonance}). 
In Figure~\ref{Fig:theta2-kappa2}, we plot values of $\rfac_2$ and $\kappa_2$ from~\eqref{r:theta2} to numerically determine the direction of bifurcation.
The fact that $\rfac_2<0$ indicates that $\rfac<\rfac_0$ along the bifurcating branch,
so we expect a supercritical Hopf bifurcation with stable periodic solutions. 
The values of $\kappa_2$ are also negative, which shows that wavenumber decreases as amplitude increases along the branch.

\subsection{The case \texorpdfstring{$\beta=r>0, \nu=0$}{β>0, ν=0}}\label{Ss.nu0}
First, we claim that no bifurcation occurs for any $\rfac>0$ if $\nu=0$ and $\beta = r\in(0,1]$.  Note
\[
H_{\beta,0}(t) := \re G_{\beta,0}(it) = \frac12\int_{-\infty}^\infty \exp\left(-\frac{|y|^{\beta+1}}{\beta+1}-ity\right)\,dy\,.
\]
For $\beta=1$ this is a Gaussian and it never vanishes. 
For $0<\beta<1$, the function $h_\beta(y)=\exp(-|y|^{\beta+1}/(\beta+1))$ is the characteristic function 
(Fourier transform) of a symmetric $\alpha$-stable distribution from probability theory, 
with $\alpha=\beta+1$, see \cite[Sec.~XVII.6]{Feller1968}.
Hence $H_{\beta,0}(t)$ is a positive multiple of this distribution's density, 
which can have no zeros due to its scaling invariance under convolution.
Hence \eqref{c:hopf1} cannot hold. 

For $\beta>1$, the inverse Fourier transform of $h_\beta$ is known \emph{not} to be a probability density,
and this means $H_{\beta,0}(t)$ \emph{must} cross zero for some $t$. Transverse crossings would then ensure bifurcation,
since $\im G_{\beta,0}(it) = -\int_0^\infty h_\beta(y)\sin(y)\,dy<0$ due to the monotonic decrease of $h_\beta$.

\subsubsection*{Odd $\beta>1$}
Strikingly, asymptotics suggests that \emph{infinitely many} bifurcations occur
for odd integers $\beta>1$ in particular. 
For in this case, $h_\beta$ is an entire function, and its Fourier transform $H_{\beta,0}(t)\to0$ 
as $|t|\to\infty$ at a super-exponential rate. A standard saddle-point analysis 
(details omitted) indicates 
\begin{equation}\label{a:Hbeta0}
H_{\beta,0}(t) \sim \sqrt{\frac{2\pi}\beta} t^{\frac{1-\beta}{2\beta} }
\exp\left( t^{\frac{\beta+1}\beta} c_\beta\right)
\cos\left( t^{\frac{\beta+1}\beta}s_\beta - \frac\pi 4 \frac{\beta-1}\beta \right)
\qquad\text{as $t\to\infty$},
\end{equation}
\[
\text{where}\qquad
c_\beta= \frac{\beta}{\beta+1} \cos \frac\pi 2\frac{\beta+1}\beta <0 ,
\qquad
s_\beta= \frac{\beta}{\beta+1} \sin \frac\pi 2\frac{\beta+1}\beta >0 .
\]
Furthermore $G_{\beta,0}(it)\sim \sfrac 1{it}$ as $t\to\infty$ (integrate \eqref{d:G} by parts).
Hence, at the bifurcation points that approximate the zeros of \eqref{a:Hbeta0}, 
the bifurcation parameter $\rfac$ and wave number in \eqref{d:waveno} satisfy
\[
\vth = \rfac^{\frac1{\beta+1}}G_{\beta,0}(0) \sim \frac 1{t^2} , \qquad 
\kappa \sim \frac{t^{-1}}{G_{\beta,0}(0)} .
\]
\begin{figure}
 \includegraphics[width=0.5\textwidth]{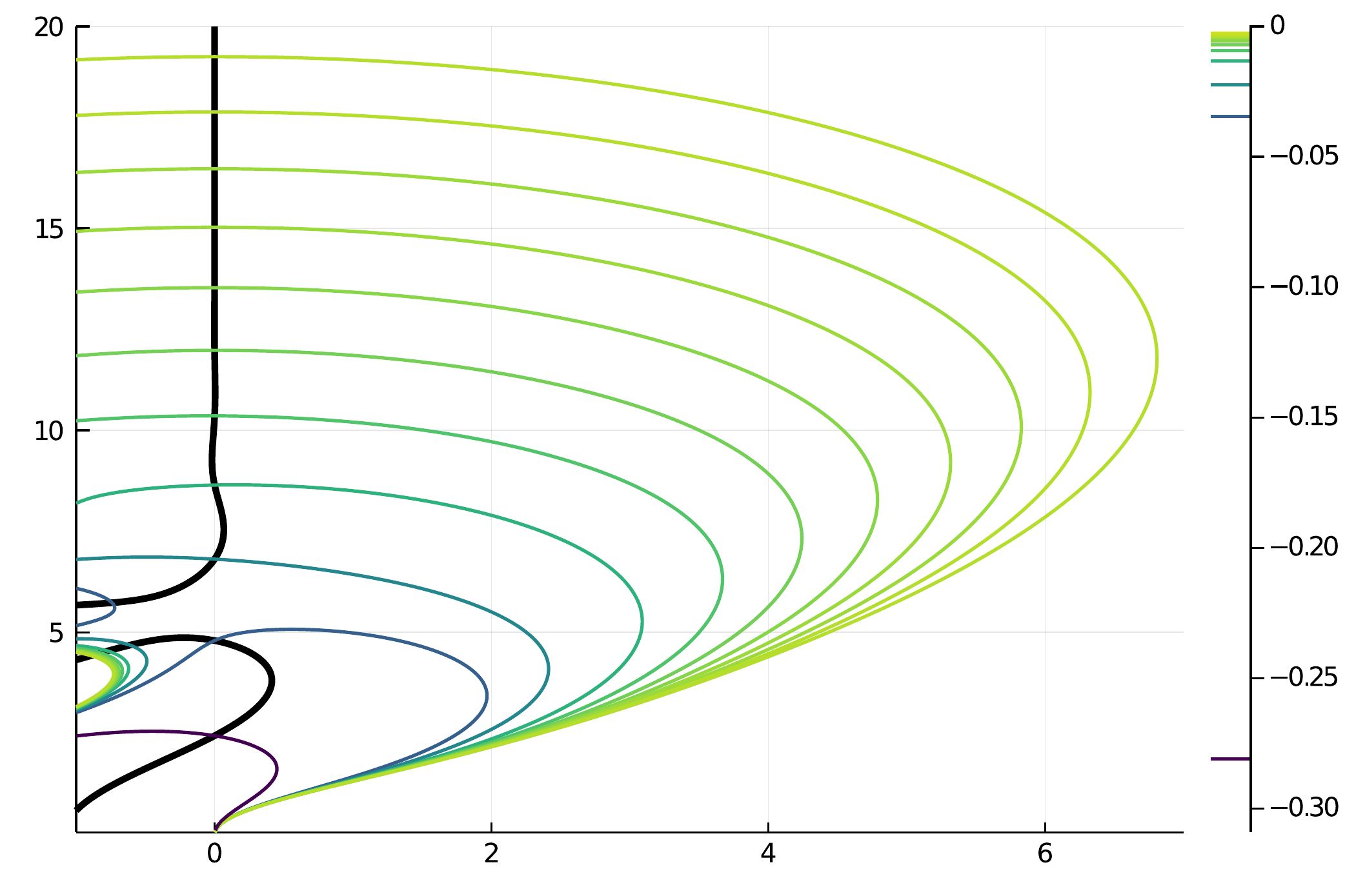}\hfill
 \includegraphics[width=0.5\textwidth]{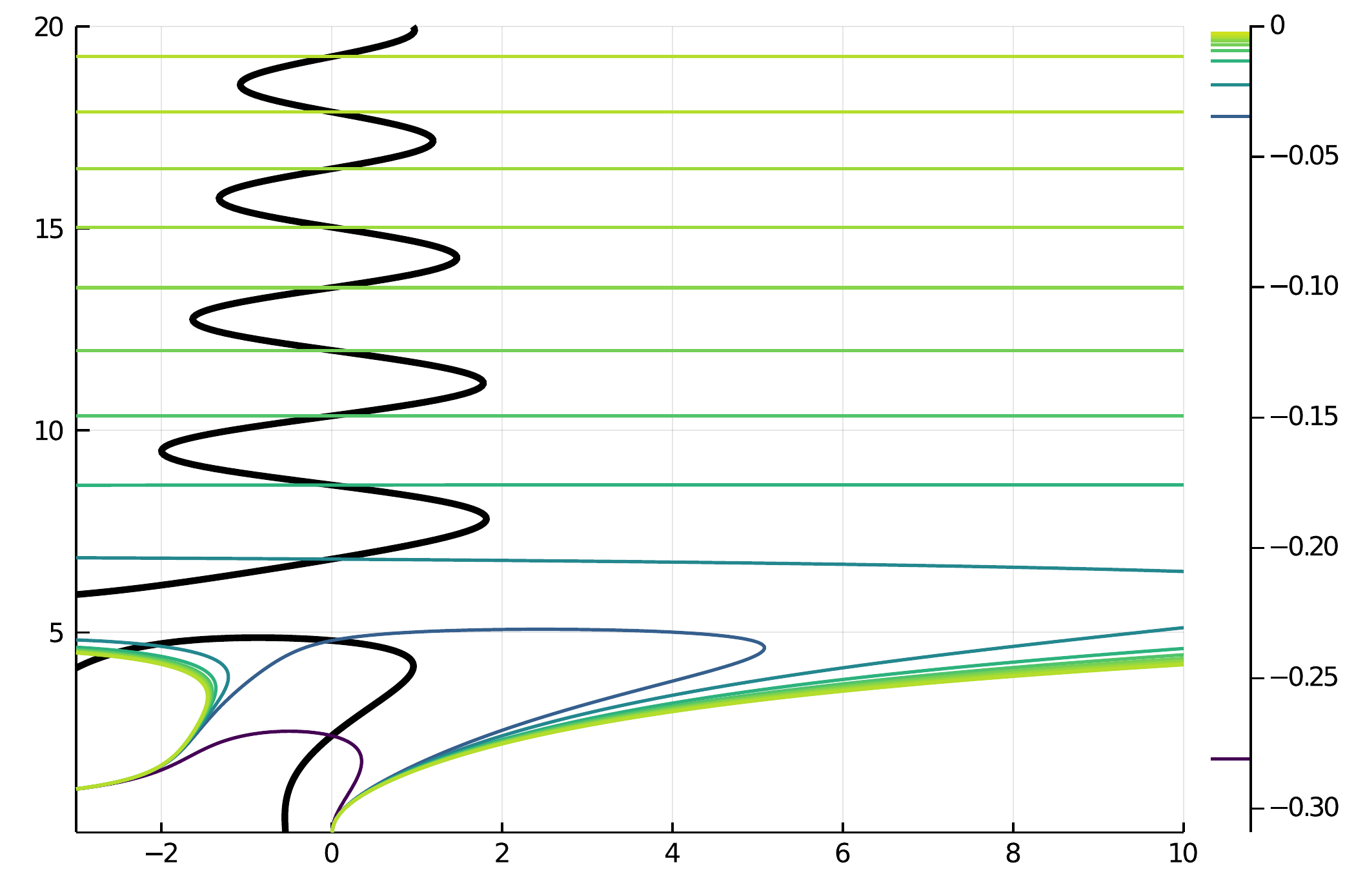}
 \caption{Contours of $\im\sfrac{G(\lambda)}{\lambda} = 0$ (black) and $\re\sfrac{G(\lambda)}{\lambda} = -\vartheta_0$ 
 for values $\vartheta_0$ given in Table~\ref{table:OddBeta3} (blue to lemon).
 (Left: unscaled. Right: real axis expanded to compensate the exponential factor in~\eqref{a:Hbeta0}) } \label{fig:contours}
\end{figure}
In Figure~\ref{fig:contours} we observe that at successive bifurcation points, evidently eigenvalues cross the imaginary axis in opposite directions,
according to \eqref{e:dlam}. Numerics\footnote{The integrals $G(\lambda)$ and roots are computed in Julia~\cite{Julia-2017} using the QuadGK~\cite{quadgk} and Roots package.} (see column \textsc{RelErr} in Table~\ref{table:OddBeta3})
indicates that for $\beta=3$, the approximation in 
\eqref{a:Hbeta0} is reasonably good already for the first zero, 
where $\re\frac{d\lambda}{d\vth}<0$ by \eqref{e:deriv:crossing}.

\begin{figure}
\pgfplotstableset{1000 sep={\,},
every head row/.style={before row=\toprule,after row=\midrule},
every last row/.style={after row=\bottomrule}}
\pgfplotstableread{BDoscillate-odd_beta3.csv}\TableOddBeta
\begin{center}
\pgfplotstabletypeset[create on use/index/.style = {create col/expr accum={\pgfmathaccuma+1}{0}},
skip rows between index={11}{20},
columns/vartheta0/.style={sci,sci zerofill,sci sep align,precision=3,sci superscript,column name=$\vartheta_0$},
columns/theta0/.style={sci,sci zerofill,sci sep align,precision=3,sci superscript,column name=$\rfac_0$},
columns/theta2/.style={sci,sci zerofill,sci sep align,precision=3,sci superscript,column name=$\rfac_2$},
columns/kappa0/.style={sci,sci zerofill,sci sep align,precision=3,sci superscript,column name=$\kappa_0$},
columns/kappa2/.style={sci,sci zerofill,sci sep align,precision=3,sci superscript,column name=$\kappa_2$},
columns/rehatg11/.style={sci,sci zerofill,sci sep align,precision=3,sci superscript,column name=$\re \hat g_1(1)$},
columns/abshatL2/.style={sci,sci zerofill,sci sep align,precision=3,sci superscript,column name=$\lvert\hat L(2)\rvert$},
columns/relerr/.style={sci,sci zerofill,sci sep align,precision=3,sci superscript,column name=\textsc{RelErr}},
columns={index,vartheta0,theta0,kappa0,theta2,kappa2,rehatg11,abshatL2,relerr}
]{\TableOddBeta}
\captionsetup{type=table}
\caption{Numerical values of the first 11 roots of $H_{3,0}(t)$, the transversal crossing condition~\eqref{e:deriv:crossing}, the nonresonance condition~\eqref{e:nonresonance} and relative error to approximate solution from~\eqref{a:Hbeta0}.}\label{table:OddBeta3}
\end{center}
\end{figure}

Thus as $\rfac$ decreases through the largest bifurcation point, one pair of eigenvalues emerges
into the right half plane and the stable constant state becomes unstable.
Decreasing across the next bifurcation point, some pair of eigenvalues must cross back into the left half plane,
and it must be the same pair. So the constant state alternates stability between successive 
bifurcation points.  

For small $\rfac$ the instabilities are extremely weak, however. 
The unstable eigenvalue becomes super-exponentially close to the imaginary axis 
since $|G'(\lambda)|$ decays only algebraically. 
Thus for small $\rfac$ the constant state becomes essentially neutrally stable.

Since Table~\ref{table:OddBeta3} indicates $\rfac_2<0$ in \eqref{r:theta2} in every case,
our study of the direction of bifurcation in section~\ref{Ss.supercritical} 
indicates that a branch of periodic solutions emerges for $\rfac$ \emph{below} each bifurcation point $\rfac_0$. 
We expect that the emerging periodic solutions alternate between stable and unstable, 
corresponding to the stability of the constant flux solution just \emph{above} $\rfac_0$.

Regarding what may happen to bifurcating solutions in the large we have no firm information.
One relatively simple possibility is that, as $\rfac$ decreases, each bifurcating branch reconnects with the next branch,
but only after folding over at some lower value of $\rfac$ and changing its stability.

\subsubsection*{Other \texorpdfstring{$\beta> 0$}{β>0} not odd}
It is possible to show that some eigenvalue crossing must always occur 
when $\beta\in(4k+1,4k+3)$ for some $k\in\N_0$.  
For in this case we can study the asymptotics of large roots of equation~\eqref{e:eval2}, 
and show the constant flux solution becomes unstable as $\rfac\rightarrow0$. 
From Section~\ref{Ss.steadyfirst} we know it is stable for large $\rfac$,
so some eigenvalue must cross the imaginary axis.

We approximate $G(\lambda)=G_{\beta,0}(\lambda)$ as $|\lambda| \rightarrow\infty$ with $\operatorname{Re}(\lambda)  \geq0$ by integrating by parts
to get
\begin{equation}\label{e:G:IntByParts}
G(\lambda)  =\frac{1}{\lambda}-\frac{1}{\lambda}\int_{0}^{\infty}
\exp\biggl(  -\frac{y^{\beta+1}}{\beta+1}\biggr)  e^{-\lambda y}y^{\beta} \, dy = \frac{1}{\lambda}+O\biggl(  \frac{1}{\left\vert
\lambda\right\vert ^{2+\beta}}\biggr) \,.
\end{equation}
Thanks to the first exponential factor in the integrand,  this expansion is valid for any $\beta\geq0$ in the region $\operatorname{Re}(\lambda)\geq-1+\delta$ for any $\delta>0$.
 We obtain the leading order approximation 
\begin{equation}\label{e:zeta:root}
\vartheta \lambda+\frac{1}{\lambda} = 0,  \qquad\text{hence}\qquad  \lambda \simeq \pm \frac{i}{\sqrt\vartheta} \qquad\text{as}\quad \vartheta\rightarrow0 \,.
\end{equation}
In order to check the stability of these roots, we need to compute the next order of $G(\lambda)$ in~\eqref{e:G:IntByParts}. 
By a change of variables $z=\lambda y$ and using  a contour deformation for large values 
of $z$ to keep the contour in the region with $\operatorname{Re}(z)\geq0$, we get 
\begin{align*}
\int_{0}^{\infty}\exp\biggl(  -\frac{y^{\beta+1}}{\beta+1}\biggr)  e^{-\lambda y} y^{\beta}dy
&=\frac{1}{ \lambda^{1+\beta}}\int_{\frac{\lambda}{\vert
\lambda\vert}\mathbb{R}}\exp\biggl( -\frac{z^{\beta+1}}{( \beta+1) \lambda  ^{1+\beta}}\biggr)  e^{-z}z^{\beta}dz \\
&\simeq\frac{1}{\lambda^{1+\beta}}\int_{0}^{\infty}e^{-z} z^{\beta}dz
=\frac{\Gamma(\beta+1)  }{\lambda^{1+\beta}}
\qquad\text{ as }\left\vert \lambda\right\vert \rightarrow\infty\,.
\end{align*}
This is  valid for any $\beta\geq 0$ in the region $\operatorname{Re}(\lambda) \geq
-1+\delta$ for $\delta>0$. 
Hence from~\eqref{e:G:IntByParts}, we get the asymptotics $G(\lambda) \simeq \lambda^{-1} -\Gamma(\beta+1) \lambda^{-(\beta+2)}$ for $|\lambda| \gg 1$. 
Here, the branch of the analytic function is chosen by extending analytically 
such that $1^{2+\beta}=1$.
We then approximate  equation \eqref{e:eval2} as
\begin{equation}\label{e:zeta:root:expansion}
\vartheta\lambda^{2}+1-\frac{\Gamma(\beta+1)}{\lambda^{1+\beta}} \simeq 0 \,.
\end{equation}
We consider  the approximation of the root with positive imaginary part in~\eqref{e:zeta:root}.
By dividing~\eqref{e:zeta:root:expansion} by $\vartheta$ and taking the square root, we get 
\begin{align*}
\lambda &= \frac{i}{\sqrt{\vartheta}} \left(  1+\frac{\Gamma(\beta+1)}{\lambda^{1+\beta}
}\right)^{\frac{1}{2}} \simeq\frac{i}{\sqrt{\vartheta}}\left(  1-\frac{\Gamma(\beta+1)}{2\lambda^{1+\beta}}\right) \qquad\text{for } |\lambda| \gg 1 \,.
\end{align*}
With the leading approximation of $\lambda$ from~\eqref{e:zeta:root}, we obtain
\begin{equation}\label{e:biglambda}
\lambda\simeq \frac{i}{\sqrt{\vartheta}} 
\left( 1-\frac{\Gamma(\beta+1) \, \vartheta^{\frac{1+\beta}{2}}}{2} \exp\biggl(-\frac{i\pi(\beta+1)}{2} \biggr) \right)
\qquad\text{as}\quad \vartheta\rightarrow 0 \,,
\end{equation}
where we used $i^{1+\beta}=e^{\frac{i\pi(\beta+1)}{2}}$. 
For $\beta=0$, we obtain $\re(\lambda)\simeq -\frac12$,
consistent with the result from Section~\ref{Ss.beta0} that the constant flux state is linearly stable. 
Regarding the general case, we find 
\begin{equation}\label{e:n:zeta:stability}
\sign\bigl(\operatorname{Re}(\lambda)\bigr)  \simeq -\sign\Biggl(\sin\biggl(\frac{\pi(\beta+1)}{2}\biggr)\Biggr) \qquad\text{ as } \vartheta \to 0 \,.
\end{equation}
Thus we obtain linear instability for $\beta\in (4k+1,4k+3)$ with $k\in \N_0$, as mentioned above. 
In case $\beta\in (4k-1,4k+1)$ for $k\in \N_0$, we have $\re(\lambda)<0$ for a root of \eqref{e:eval2} satisfying \eqref{e:biglambda}, but we do not know about other roots.
Nevertheless, this analysis suggests that in the limit $\vartheta\to0$, the constant flux state changes its stability as $\beta$ passes an odd integer, which is exactly when the system shows infinitely many bifurcations as discussed in Section~\ref{Ss.nu0}.

\subsubsection*{The case \texorpdfstring{$\beta\geq 0$}{β≥0} not odd and \texorpdfstring{$\nu\ll 1$}{ν≪1}}

We close this section by indicating that an analysis is possible for the case $\beta\geq 0$ not odd with $\nu\ll 1$ but positive. 
In \AppendixSupp~\ref{Ss.nu_small}, 
we arrive, consistent with~\eqref{e:zeta:root}, at the approximation of the root $\lambda$ of~\eqref{e:eval2} with positive imaginary real part to second order as
\[
  \lambda\simeq\left(  \frac{\Gamma(\nu+1)}{\vartheta}\right)^{\frac{1}{\nu+2}}  \, (i + \lambda_2) \qquad \text{as }\ \vartheta\rightarrow 0 ,
\]
with $\lambda_2 = O(\vartheta^{\frac{\beta+1}{\nu+2}})$. The real part of $\lambda_2$ determines the stability of the stationary solution. Approximating $\lambda_2$ in the limit $\nu\to 0$, we find
\begin{equation}\label{e:null1:zeta2}
  \operatorname{Re}(\lambda_2)  \simeq \frac{\pi \nu}{4}-\frac{\Gamma(\beta+1)  }{2
}\vartheta^{\frac{\beta+1}{\nu+2}}\sin\Bigl(
\frac{\pi (  \beta+1)  }{2}\Bigr) \qquad\text{as }\ \nu\to 0 .
\end{equation}
First, for $\nu=0$, we recover the stability condition~\eqref{e:n:zeta:stability}. Next, for $\nu>0$, 
the stationary solution becomes unstable as soon as $\vartheta$ is small enough and oscillatory behavior can be expected, in accordance to the asymptotics from Section~\ref{Ss.supercritical}. Finally, we obtain a curve of stability with $\nu$ depending on $\vartheta$ if $\sin\bigl(\frac{\pi(  \beta+1)}{2}\bigr) >0$, that is $\beta\in (4k-1,4k+1)$ for some $k\in \N_0$.

% Appendices were here

\appendix

\section{Moment models}\label{Ss.momentmodels} 
% \footnote{\BN{This section is new to me and I do not really understand its purpose and also not all details.}}
% 
% In the simplified model we used a change of variables to obtain transport at
% constant speed. The final model, after the rescalings, but written in the
% original variables (i.e. $y$ is the particle volume) is:%
% \begin{align*}
% \partial_{t}f\left(  y,t\right)  +\partial_{y}\left(  y^{\alpha}f\left(
% y,t\right)  \right)    & =-ry^{b}f\left(  y,t\right)  \\
% f\left(  y,t\right)    & \simeq\frac{e^{u}}{y^{\alpha}}\ \ \text{as\ \ }%
% y\rightarrow0^{+}\\
% \partial_{t}u  & =1-\int_{0}^{\infty}y^{\alpha}f\left(  y,t\right)  dy
% \end{align*}
% 
In \cite[p. 293ff]{Friedlander2000} a model is suggested for the distribution of diameter sizes of supercritical particles in the presence of homogeneous and heterogeneous nucleation having similarities to the limit model in Section~\ref{Ss.LimitModel}.
% In~\cite[p. 293ff]{Friedlander2000} a moment model is derived \footnote{\BN{Derived from what?}} based on the density in terms of the diameter, instead of the volume as it is in the case of this paper.
Homogeneous nucleation is modeled using an Arrhenius law, and specific power laws for the particle growth and the particle removal are assumed. 
In~\cite[(10.46)]{Friedlander2000} constant 
growth of this diameter is assumed, which translates for our coefficients~\eqref{coeffassum} to the case $\alpha=\sfrac{2}{3}$. Part of the analysis of this case can also be found in~\cite{PratsinisFriedlanderPearlstein1986}.
% This density would satisfy:%
% \[
% n\left(  d_{p}\right)  d\left(  d_{p}\right)  =f\left(  y\right)  dy
% \]
% 
% 
% The rate of growth used in Friedlander (see (10.46)) is:%
% \[
% \frac{d\left(  d_{p}\right)  }{dt}=1
% \]
% 
% 
% Strictly speaking the assumption is that the rate of change of the particle
% diameter is proportional to the supersaturation $\left(  S-1\right)  ,$ but in
% our model this is almost constant (in the limit under consideration). Then:%
% \[
% \frac{dy}{dt}=y^{\frac{2}{3}}%
% \]
% 
% 
% Then:%
% \[
% \alpha=\frac{2}{3}%
% \]

% 
% Notice that in the limit under consideration the rate of growth is independent
% of the supersaturation, because this is almost constant.
% 
% \bigskip
For this choice of the coefficients, the limit model~\eqref{eq:limit:u}--\eqref{eq:limit:boundary} becomes
\begin{align}
\partial_{t}f(x,t)  +\partial_{x}\bigl(  x^{\frac{2}{3}}f(
x,t)\bigr)    & =-\rfac x^{\rexp}f(  x,t)\,,  \label{moment:transportf}\\
x^{\frac{2}{3}}f(x,t)    &\simeq e^{u}\ \ \text{as\ \ }%
x\rightarrow0^{+}\,. \label{moment:boundary}
\end{align}%
Following~\cite{PratsinisFriedlanderPearlstein1986}, the authors use three moments, denoted as $N,\ A,\ R$ which are the number of clusters, the area and the 
radius, respectively. In
our notation these moments are
\[
N=\int_{0}^{\infty}f(x,t)  dx\,,\quad A=\int_{0}^{\infty}%
x^{\frac{2}{3}}f(x,t)  dx\,,\quad  R =\int_{0}^{\infty}%
x^{\frac{1}{3}}f(x,t)  dx \, .
\]
It is then possible to calculate the evolution equations for these moments,
% \begin{align*}
% \partial_{t}N  & =e^{u}-r\int_{0}^{\infty}x^{b}f\left(  x,t\right)  dx\\
% \partial_{t}A  & =\frac{2}{3}\int_{0}^{\infty}x^{\frac{1}{3}}f\left(
% x,t\right)  dx-r\int_{0}^{\infty}x^{b+\frac{2}{3}}f\left(  x,t\right)
% dx=\frac{2}{3}M_{1}-r\int_{0}^{\infty}x^{b+\frac{2}{3}}f\left(  x,t\right)
% dx\\
% \partial_{t}M_{1}  & =\frac{1}{3}N-r\int_{0}^{\infty}x^{b+\frac{1}{3}}f\left(
% x,t\right)  dx
% \end{align*}
which we have to complement with the equation for the monomer
concentration
\[
\partial_{t}u=1-\int_{0}^{\infty}x^{\frac{2}{3}}f(x,t)  dx=1-A .
\]
In order to obtain a closed system of ODEs we need to make the removal term precise. In~\cite{Friedlander2000}, it is assumed that $\rfac=0$, which gives the following system for the moments:
\[
\partial_{t}N=e^{u}\,,\qquad  \partial_{t}A=\frac{2}{3}R\,, \qquad  \partial_{t}R=\frac{1}{3}N\,, \qquad  \partial_{t}u=1-A .
\]
This is a system which is not able to generate oscillations, since the moments just grow in time. 

In~\cite{PratsinisFriedlanderPearlstein1986} a removal mechanism is suggested that eliminates 
clusters with a mean life time  $\tau>0$, which is equivalent to choosing $\rexp=0$ and $\rfac=\frac{1}{\tau}>0$ in~\eqref{moment:transportf}.
With this choice, the ODEs for the moments take the form
% One of the equations (for $A$) and this removal rate is mentioned in the book of Friedlander, but
% the details about the equation are in the paper. 
\[
\partial_{t}N=e^{u}-\rfac N\,, \qquad \partial_{t}A=\frac{2}{3}R-\rfac A\,, \qquad  \partial_{t}R=\frac{1}{3}N-\rfac R\,, \qquad  \partial_{t}u=1-A \,.
\]
This ODE system derived from~\eqref{moment:transportf} and~\eqref{moment:boundary} is almost the same as the system considered in~\cite{PratsinisFriedlanderPearlstein1986} with only two minor differences: First, in \cite{PratsinisFriedlanderPearlstein1986} there is a term associated to the flux of area and radius at the critical radius present. In our model a similar term is shown to be negligible in the limit that we consider (see~\eqref{e:negligible_fluxes}). Second, for the flux of clusters, as denoted in~\cite{Friedlander2000} by $I$, the full Arrhenius formula is used, whereas we obtain the exponential approximation~\eqref{Jinftydef} leading to the boundary condition~\eqref{moment:boundary}.

% end of appendix on moments

%%% Goes to supplementary material 

%\appendix
% \section{Appendix}

\section{Time scale for the quasistationary approximation}\label{sec:quasistationary}
%\marginnote{\RP{Can this go to supplementary material?}}
%
In this section we are going to argue that solutions $\{  n_{k}\}$ to the Becker-D\"oring equations~\eqref{n1eq}, \eqref{nkeq} with fluxes~\eqref{jkeq} and coefficients~\eqref{coeffassum} stabilize to a quasistationary distribution for
values of $k$ of the order of the critical size. Given that the critical
size rescales like $k_{\crit}=(q/\eps) ^{\frac{1}{\gamma}}$ we expect this stabilization time to be an algebraic function of
$\eps^{-1}$.

We expect that  the concentration of monomers $n_{1}$ changes on a time scale  which is exponentially large in
$\eps$  and that  $n_{1} \simeq \bar n_1=1+\eps \simeq 1$ during  the whole evolution. Since the removal term $R$ becomes in the variable $\rfac= TX^\rexp R$ and the scaling~\eqref{eq:scales} from Section~\ref{Ss.scales} exponentially small in $\eps$, 
we need to estimate the characteristic time scale for the stabilization to the stationary solution of
\begin{align}
\partial_{t}n_{k}  &  =(a_{k-1}\bar{n}_{1}n_{k-1}-b_{k}n_{k})-(a_{k}\bar{n}_{1}n_{k}-b_{k+1}n_{k+1})\,, \qquad  k\geq
2\label{A1E}
%\\
% {\bar  n}_{1}&=1 \,. \label{A2E}%  % RP: Either 1+\eps, or just leave out, I'd say
\end{align}
Given that the
%solutions of \eqref{A1E}, \eqref{A2E} satisfy a comparison principle, it is
solutions of \eqref{A1E} satisfy a comparison principle, it is
enough to estimate the time scale of convergence to the steady state of the solutions for initial values that are a
constant, either positive or zero. Using \eqref{coeffassum},
 approximating  $n_{k}$ by a differentiable function $n(k,t)$ and  replacing differences by continuous derivatives
we obtain an initial-boundary value problem,
\begin{align}
\partial_{t}n(k,t)  &=-\partial_{k}\bigl( k^{\alpha}(  \eps-q k^{-\gamma})
n(k,t)  \bigr)  +\partial_{k}^{2}\bigl(  k^{\alpha}(  1+q k^{-\gamma})  n(k,t)  \bigr) \,, \quad k \geq 1\,, \label{A8} \\
n(  1,t) & =1\,. \label{A8a}%
\end{align}
This approximation will not  yield an accurate quantitative approximation of  
%solutions of \eqref{A1E}, \eqref{A2E}, 
solutions of \eqref{A1E}, 
but it suffices to estimate the characteristic time scale for stabilization to a steady state. 
(See~\cite{Velazquez1998,ConlonSchlichting2019} for similar approximations of the Becker--D{\"o}ring system.)
For that we must decompose the process into several stages. 

For initial data  $n(k,0)=0$, the first stage in 
the build-up of the stationary solution is the stabilization of the region
where $k$ is of order one. If $k$ is bounded we can approximate \eqref{A8} by
\begin{equation}
\partial_{t}n(k,t)  =\partial_{k}\bigl(  q k^{\alpha-\gamma}n(k,t)  \bigr)  +\partial_{k}^{2}\bigl( k^{\alpha}(  1+q k^{-\gamma}) n(  k,t)  \bigr)\,, \qquad k > 1\,. \label{A8b}%
\end{equation}
           Solutions of \eqref{A8b}, \eqref{A8a} with zero initial data stabilize,
           satisfying
\begin{align}
n(k,t) &\simeq n_s(k) \qquad \text{ for } t \gg 1\,\label{A8d}
\end{align}
with 
\begin{align}\label{A8e}
n_{s}( k) &= \frac{\exp\Bigl( - \int_1^k \frac{q}{\xi^{\gamma} + q}d\xi\Bigr)}{k^{\alpha} \big( 1 + \frac{q}{k^{\gamma}}\big)} \simeq\frac{1}{ k^{\alpha}}\exp\Bigl( - \frac{q}{1{-}\gamma}
k^{1{-}\gamma}\bigl(1+ o(1)\bigr) \Bigr)      \quad \text{ for } k \gg 1\,.
\end{align}

If $n(k,0)  =C>0$ there is a transient evolution (algebraically slow in $\eps$) before we  obtain  \eqref{A8d}. To describe this we recall $k_\crit$ from~\eqref{criticalsize} and introduce
\begin{equation}
k=k_\crit \, x\,, \qquad  n(k,t)  =f(x,s) \,, \qquad  t=\frac
{1}{\eps}\Big(  \frac{\eps}{q}\Big)  ^{\frac{\alpha
}{\gamma}-1}s\,. \label{A9}%
\end{equation}
With $\delta =q^{-\frac{1}{\gamma}}  \eps^{\frac{1-\gamma}{\gamma}}$ equation \eqref{A8} 
becomes 
\begin{equation}
\partial_{s}f(x,s)  =\delta\partial_{x}^{2}\bigl(  x^{\alpha}(  1+\eps x^{-\gamma})  f(x,s) \big)  -\partial_{x}\big(  x^{\alpha}(  1-x^{-\gamma})  f(x,s)  \bigr)\,, \quad  x > k_\crit^{-1}\,,\label{A7}%
\end{equation}
with boundary conditions, due to \eqref{A8a},  
\begin{equation}
f(k_\crit^{-1},s)  =1\,, \qquad s>0\,, \label{A7a}%
\end{equation}
and  initial conditions
\begin{equation}
f(  x,0)  =C\,, \qquad x>k_\crit^{-1}\,. \label{A7b}%
\end{equation}
For $x \in (0,1)$ 
we can approximate \eqref{A7}--\eqref{A7b}  
by
\begin{equation}
\partial_{s}f(x,s)  =-\partial_{x}\bigl(x^{\alpha}(  1-x^{-\gamma})  f(x,s)  \bigr)\,.  \label{A7c}%
\end{equation}
                    For \eqref{A7c} the term $x^{\alpha}\big(1-\frac{1}{x^{\gamma}}\big)  f(x,s)  $ is
constant along characteristics. That is, 
suppose that we
denote by $X(\xi,s;x)  $ the unique solution of
\[
\frac{\partial X}{\partial\xi}= X^{\alpha}\big(  1- X^{-\gamma}\big )  \qquad \text{ for  } 0\leq\xi\leq
s\,, \qquad  X(s,s;x)  =x \in (0,1)\,,
\]
                                then
\[
x^{\alpha}(  1-x^{-\gamma})  f(x,s)  =  X(0,s;x)^{\alpha}\bigl(  1-X(0,s;x)^{-\gamma}\bigr) C\,.
\]
It is easily seen that $X(0,s;x) \to 1$ as $s \to \infty$ and therefore
\begin{equation}
\lim_{s\rightarrow\infty} \sup_{x\in(0,b)  }\big( x^{\alpha-\gamma}f(x,s)\big  ) =0  \qquad \text{ for any } b \in (0,1)\,.
\label{A7e}%
\end{equation}
Formula \eqref{A7e} holds for the solutions of \eqref{A7c}. In the case of the
solutions of the parabolic equation \eqref{A7} we can only expect to have for some $a>0$
\begin{equation}
 x^{\alpha-\gamma}f(x,s)  =O\big(e^{-as}\big) \qquad \text{ as } s\rightarrow\infty \quad \text{ for } k\gg
1\,, \;  x\ll1 \,.\label{A7f}%
\end{equation}
Using \eqref{A9} it follows that in times $t$ algebraic in
$\eps^{-1}$ we have that $n(k,t)$ is small if $k$ is large.
In particular notice that $n(k,t)  $ is small if $s\gg1$ and $x$
is of order one, i.e. $k$ of order $\eps^{-\sfrac{1}{\gamma}}$.

In order to describe $n(k,t)  $ for  $k$ of order one we
use \eqref{A8b}, which must be solved, due to \eqref{A7f} assuming that
$k^{\alpha-\gamma}n(k,t)  \rightarrow0$ for
$k\gg1$ (and $x\ll1$) and $t$ sufficiently large (algebraic), say of order
$\eps^{\frac{\alpha}{\gamma}-2}$. The solutions of
\eqref{A8a}, \eqref{A8b} with this condition approach, after the first
transient state in an additional time $t$ of order one, to $n_{s}(k)  $ in \eqref{A8e}. In particular, \eqref{A8d}, \eqref{A8e} hold after
this transient regime.

We now describe how $n(k,t)  $ evolves towards a  a steady state for $k\simeq k_\crit$. We use the WKB method and define 
 $S(x,s)$ via 
\begin{equation}
f(x,s)  =\exp\bigl(  \delta^{-1} S(x,s)\bigr)\,.  \label{B3}%
\end{equation}
Plugging \eqref{B3} into \eqref{A7} we obtain to  leading order 
 the  Hamilton-Jacobi equation 
\begin{equation}
\frac{\partial S}{\partial s}=x^{\alpha}\Big(
\frac{\partial S}{\partial x}\Big)  ^{2}-x^{\alpha}\Big(1-\frac{1}{x^{\gamma}}\Big)  \frac{\partial S}{\partial x}\,.
\label{B4}%
\end{equation}
The asymptotic formulas \eqref{A8d}, \eqref{A8e} yield the  matching
condition
\begin{equation}
S(x,s)  \simeq -\frac{x^{1-\gamma}}{1{-}\gamma}\,, \qquad \text{ as } x\rightarrow0^{+}\,. \label{B4a}%
\end{equation}
Equation \eqref{B4} defines, with $p=\frac{\partial S}{\partial x}$, the  Hamiltonian system
\begin{align}
\frac{dx}{ds}  & =\frac{\partial
H}{\partial p}\,, \qquad 
\frac{dp}{ds}    =
-\frac{\partial H}{\partial x} \label{H2}%
\qquad \text{ with } \qquad H=-x^{\alpha}p^{2}+x^{\alpha}\big(1-x^{-\gamma}\big)   p\,.
                    \end{align}
Then, equation \eqref{B4} yields that along characteristics we have
\begin{equation}
\frac{dS}{ds}=-x^{\alpha}p^{2}\,. \label{H3}%
\end{equation}
                   We normalize by subtracting  from $s$ the time required to build the
asymptotics \eqref{A8d} and  obtain from  \eqref{B4a} that
\begin{equation}
p\simeq -x^{-\gamma}\qquad \text{ as }s\rightarrow0 \text{ and }  x\rightarrow0\,. \label{H4}%
\end{equation}
 From \eqref{B4a} we deduce that $\frac{\partial S}{\partial s}=0$ as $x \to 0$ and hence             the trajectories lie in the level set
 $\{H=0\}$.  Consequently,
the characteristics starting at the line $x=0$ satisfy 
$-x^{\alpha}p^{2}+x^{\alpha}\big(1-x^{-\gamma}\big)   p=0$  such that 
along characteristics $p=(1-x^{-\gamma})$ and
 \[
\frac{dx}{ds}=  x^{\alpha}\big( 1- x^{-\gamma}\big) - 2x^{\alpha} p = -x^{\alpha}\big(1-x^{-\gamma}\big) = - x^{\alpha} p\,.
\]
Using \eqref{H3} we obtain
$\frac{dS}{dx}=\frac{x^{\alpha}p^{2}}{x^{\alpha}p}=p=1-x^{-\gamma}$
and as a consequence 
\begin{equation}
S(x,\infty)  =-\frac{x^{1-\gamma}}{1{-}\gamma}+x\,, \qquad  0<x<1\,. \label{H5}%
\end{equation}
The characteristics approach asymptotically to $(x,p)  =(1,0)  $ as $s\rightarrow\infty$.

It remains to estimate the time scale for the solution to overcome the point
$x=1$. The characteristics arrive to  $x=1$ in times $s$ of order
one, which means in algebraic times in the original variable $t$.
In the region $x\simeq 1$ we need to use again the full diffusion equation. 
%The characteristics $x(s)  $ approaches exponentially to $x=1$ as $s\rightarrow\infty$.  
In order to understand the behavior of $n(k,t)$ near $k_{\crit}$ we introduce a new
variable $y= k_{\crit}^{-\frac{\gamma+1}{2}}(k-k_{\crit})$. The characteristics   arrive to some $k=k_{\crit}%
-L(k_{\crit})^{\frac{\gamma+1}{2}}$ with $L$ large in
times bounded
algebraically in the variable $t$. It is important to take into account that
$  k_{\crit}^{\frac{\gamma+1}{2}}\ll k_{\crit},$ because due to
this, the region where we use this approximation is a small layer compared
with the region described by the critical size. The Hamilton-Jacobi
approximation is valid as long as $\vert y\vert \gg1$ (with $y<0$).

Linearizing in \eqref{A8} and writing $N(y,t)  =n(k,t)  $ we obtain the approximation
\begin{equation}
\partial_{t}N(y,t)  =  k_{\crit}^{\alpha-\gamma
-1}\Bigl(  \partial_{y}^{2} N(y,t)   -\gamma
q\partial_{y}\bigl(  yN(y,t)  \bigr) \Bigr)  \label{H6}%
\end{equation}
Using \eqref{B3}, \eqref{H5} we obtain 
\[f(x,s)   =\exp\Bigl(  \frac{1}{\delta}S(x,s)\Bigr) \simeq e^{-\frac{\gamma}{(1-\gamma)\delta}} \exp\Bigl(  \frac{\gamma}{2\delta}(x{-}1)^{2}\Bigr)\qquad \text{ for } s\gg 1
\]
and the following asymptotics 
 in the region $y<0,$ $\vert y\vert \gg1:$%
\begin{align}
N(y,t) =f(x,s) & \simeq e^{-\frac{\gamma}{(1{-}\gamma)\delta}}   \exp\Bigl(  \frac{\gamma}{2\delta
}(x-1)  ^{2}\Bigr)=e^{-\frac{\gamma}{(1{-}\gamma)\delta}}   \exp\Bigl(  \frac{\gamma q}{2}y^{2}\Bigr)\,.  \label{H7}%
\end{align}
We can expect that the solution to \eqref{H6}, \eqref{H7} stabilizes to a steady state on a time scale that is still algebraic in $\eps^{-1}$.
Given that the initial datum is bounded, it is natural to
expect that the steady state stabilizes to a bounded steady state.
That is we solve 
 $\partial_{y} N_{s}(y)  -\gamma qyN_{s}(y)  =-J$
which gives
\begin{equation}\label{H7b}
N_{s}(y)  =J\exp\Bigl(  \frac{\gamma q}{2}y^{2}\Bigr)
\int_{y}^{\infty}\exp\Bigl(  -\frac{\gamma q}{2}\xi^{2}\Bigr)  d\xi \,.
\end{equation}
We have
\[
N_{s}(y)  \simeq J\exp\Bigl(  \frac{\gamma q}{2}y^{2}\Bigr)
\int_{-\infty}^{\infty}\exp\Bigl(  -\frac{\gamma q}{2}\xi^{2}\Bigr)
d\xi
 \simeq\sqrt{\frac{2\pi}{\gamma q}}J\exp\Bigl(
\frac{\gamma q}{2}y^{2}\Bigr)  \quad \text{ as } y\rightarrow-\infty\,.
\]
Then, the matching condition \eqref{H7} implies $J=\sqrt{\frac{\gamma q}{2\pi}}    e^{-\frac{\gamma}{(1{-}\gamma)\delta}}$. 

Finally, we can now obtain the propagation of the solution towards the region
$y\gg1$ (and also $x>1$).
We have, due to \eqref{H7b}, that $N_{s}(y)  \simeq\frac{J}{\gamma qy}$ as 
$y\rightarrow\infty$. 
In this region ($x>1$) due to the absence of exponential terms in the
solution, we can transport the solution by characteristics, using \eqref{A7c}
for $x>1$. The problem describing the evolution of the concentration of
clusters with $x>1$ is:%
\begin{align*}
\partial_{s}f(x,s)    & =-\partial_{x}\big(x^{\alpha}\big(  1-x^{-\gamma}\big)  f(x,s)  \big) \,, \qquad  x>1\,,\\
f(x,s)    & \simeq\frac{J}{\gamma q}\Big(  \frac{\eps%
}{q}\Big)  ^{\frac{1-\gamma}{2\gamma}}\frac{1}{
(x{-}1)}\qquad \text{ as } x\rightarrow1^{+}\,.
\end{align*}
The speed of the characteristics is of order one in $s$, hence the stabilization to the equilibrium distribution
takes place also in algebraic times in the variable $t$.

 % arXiv
%
\section{Numerical verification of transversality and non-resonance conditions}\label{supp:num-ver}

\pgfplotstableset{1000 sep={\,},
	every head row/.style={before row=\toprule,after row=\midrule},
	every last row/.style={after row=\bottomrule}}

\pgfplotstableread{BDoscillate-resonance-transversality.csv}\TableResonanceTransversality

\begin{center}
	\pgfplotstabletypeset[columns/beta/.style={column type=r,column name=$\beta$},
	columns/rehatg111/.style={sci,sci zerofill,sci sep align,precision=3,sci superscript,column name=$\nu{=}\frac19$},
	columns/rehatg112/.style={sci,sci zerofill,sci sep align,precision=3,sci superscript,column name=$\nu{=}\frac14$},
	columns/rehatg113/.style={sci,sci zerofill,sci sep align,precision=3,sci superscript,column name=$\nu{=}\frac37$},
	columns/rehatg114/.style={sci,sci zerofill,sci sep align,precision=3,sci superscript,column name=$\nu{=}\frac23$},
	columns/rehatg115/.style={sci,sci zerofill,sci sep align,precision=3,sci superscript,column name=$\nu{=} 1$},
	columns/rehatg116/.style={sci,sci zerofill,sci sep align,precision=3,sci superscript,column name=$\nu{=}\frac32$},
	columns/rehatg117/.style={sci,sci zerofill,sci sep align,precision=3,sci superscript,column name=$\nu{=}\frac73$},
	columns/rehatg118/.style={sci,sci zerofill,sci sep align,precision=3,sci superscript,column name=$\nu{=}4$},
	columns/rehatg119/.style={sci,sci zerofill,sci sep align,precision=3,sci superscript,column name=$\nu{=}9$},
	columns={beta,[index]10,[index]11,[index]12,[index]13,[index]14,[index]15,[index]16,[index]17,[index]18},
	row predicate/.code={% Include only every 5th row
		\pgfplotstablegetelem{#1}{beta}\of{\TableResonanceTransversality}
		\pgfmathparse{int(Mod(multiply(10,\pgfplotsretval),5))}
		\ifnum\pgfmathresult=0\relax
		\else\pgfplotstableuserowfalse
		\fi}
	]{\TableResonanceTransversality}
	\captionsetup{type=table}
	\caption{Numerical verification of the transversality condition~\eqref{c:hopf2} with values of $\re \hat g_1(1)$.}\label{table:transversality}
\end{center}

%$\nu\in \{ \frac19 \frac14 \frac37 \frac23 1 \frac32 \frac73 4 9 \}
\begin{center}
	\pgfplotstabletypeset[columns/beta/.style={column type=r,column name=$\beta$},
	columns/abshatL21/.style={sci,sci zerofill,sci sep align,precision=3,sci superscript,column name=$\nu{=}\frac19$},
	columns/abshatL22/.style={sci,sci zerofill,sci sep align,precision=3,sci superscript,column name=$\nu{=}\frac14$},
	columns/abshatL23/.style={sci,sci zerofill,sci sep align,precision=3,sci superscript,column name=$\nu{=}\frac37$},
	columns/abshatL24/.style={sci,sci zerofill,sci sep align,precision=3,sci superscript,column name=$\nu{=}\frac23$},
	columns/abshatL25/.style={sci,sci zerofill,sci sep align,precision=3,sci superscript,column name=$\nu{=} 1$},
	columns/abshatL26/.style={sci,sci zerofill,sci sep align,precision=3,sci superscript,column name=$\nu{=}\frac32$},
	columns/abshatL27/.style={sci,sci zerofill,sci sep align,precision=3,sci superscript,column name=$\nu{=}\frac73$},
	columns/abshatL28/.style={sci,sci zerofill,sci sep align,precision=3,sci superscript,column name=$\nu{=}4$},
	columns/abshatL29/.style={sci,sci zerofill,sci sep align,precision=3,sci superscript,column name=$\nu{=}9$},
	columns={beta,[index]1,[index]2,[index]3,[index]4,[index]5,[index]6,[index]7,[index]8,[index]9},
	row predicate/.code={% Include only every 5th row
		\pgfplotstablegetelem{#1}{beta}\of{\TableResonanceTransversality}
		\pgfmathparse{int(Mod(multiply(10,\pgfplotsretval),5))}
		\ifnum\pgfmathresult=0\relax
		\else\pgfplotstableuserowfalse
		\fi}
	]{\TableResonanceTransversality}
	\captionsetup{type=table}
	\caption{Numerical verification of the nonresonance condition~\eqref{e:nonresonance} with values of $\lvert \hat L(2)\rvert$.}\label{table:nonresonance}
\end{center}
 % arXiv
%
\section{The case \texorpdfstring{$\nu\ll 1$}{ν≪1} and \texorpdfstring{$\beta\geq0$}{b≥0}}\label{Ss.nu_small}
%\marginnote{\AS{Clarify correct regime for $\beta$}}

%\marginnote{\RP{query: ``for any $\nu=0$??}}
It is interesting to examine this regime, given that for each $\nu>0$ we expect a
Hopf bifurcation but this is lost when $\nu=0$ for $\beta \in (0,1]$ as discussed in Section~\ref{Ss.nu0}. 
We examine first the roots yielding the Hopf bifurcation, but it is also of interest to
study the limits $\nu\to 0$ and also $\rfac\rightarrow0$, since this could yield
bifurcations with very large period and result in some interesting limiting behavior.

The starting point is the eigenvalue equation~\eqref{e:eval2}, which we repeat here as
\begin{equation}\label{supp:e:eval2}
	\vth\lambda + G_{\beta,\nu}(\lambda) = 0,
	\qquad\text{where}\qquad
	\vth = \rfac^{\frac1{\beta+1}}  G_{\beta,\nu}(0)\,,
\end{equation}
where we the function $G_{\beta,\nu}(\lambda)$ takes with the substitution $\xi = \lambda \, z$ the form
\begin{equation}\label{e:supp:def:G}
	G_{\beta,\nu}(\lambda)  =\int_{0}^{\infty}y^{\nu}\exp\biggl(  -\frac{z^{\beta+1}}%
	{\beta+1}\biggr)  e^{-\lambda z} dz  =\frac{1}{\lambda^{\nu+1}}\int_{0}^{\infty}\xi^{\nu}%
	\exp\biggl(  -\frac{\xi^{\beta+1}}{(  \beta+1)  \lambda^{\beta+1}}\biggr)
	e^{-\xi}d\xi .
\end{equation}
We note that $G_{\beta,\nu}(0)= (1+\beta)^{-\frac{\beta-\nu}{1+\beta}} \Gamma\bigl(\frac{1+\nu}{1+\beta}\bigr)>0$ and the specific value is not of importance in the following. 
We examine the position of the roots of~\eqref{supp:e:eval2} as $\vartheta\rightarrow0$. To this end, we compute the asymptotics of $G_{\beta,\nu}(\lambda)  $ as $\left\vert \lambda\right\vert \rightarrow\infty$ and obtain 
\begin{equation}\label{e:supp:G:asymp}
	G_{\beta,\nu}(  \lambda)  \simeq \frac{1}{\lambda^{\nu+1}}\int_{0}^{\infty}\xi
	^{\nu}e^{-\xi}d\xi=\frac{\Gamma(  \nu+1)  }{\lambda^{\nu+1}}%
	\qquad \text{as } \vert \lambda\vert \rightarrow\infty, 
\end{equation}
which hold true in the region $\operatorname{Re}(  \lambda)  >-1$. 
With~\eqref{e:supp:G:asymp}, we approximate the eigenvalue equation~\eqref{supp:e:eval2} for $\lambda$ as
\begin{equation}\label{e:supp:EV:zeta:asymp}
	\vartheta\lambda\simeq -\frac{\Gamma(  \nu+1)  }{\lambda^{\nu+1} }\qquad\text{as\ }\left\vert
	\lambda\right\vert \rightarrow\infty .
\end{equation}
Hence, we obtain the approximation $\vartheta \lambda^{\nu+2}=-\Gamma(\nu + 1)$ as $\rfac\rightarrow0$, which gives the  approximation of the roots of~\eqref{supp:e:eval2} as
\begin{equation}\label{e:supp:EV:zeta:asymp_val}
	\lambda \simeq \frac{\Gamma(\nu+1)^{\frac{1}{\nu+2}}}{\vartheta^{\frac{1}{\nu+2}}}e^{\pm\frac{i\pi}{\nu+2}} \qquad\text{as} \lambda \to \infty.
\end{equation}
For $\nu>0$, these roots have positive real part and hence the constant flux
solution becomes unstable for any $\nu>0$ fixed and any $\beta\geq 0$. On the other
hand if $\nu=0$ we have stability of this root if $\beta\in (4k-1,4k+1)$ for $k\in \N_0$ in agreement to the previous discussion in Section~\ref{Ss.nu0}.

Now, we compute the line in the plane $(\nu,\beta)$ which separates the
stability and instability regions as $\rfac\rightarrow0^{+}$. To this end we go back to~\eqref{e:supp:def:G} and expand the exponential term
\begin{equation}\label{e:supp:G:exp_expand}
	G_{\beta,\nu}(  \lambda)  =\frac{1}{\lambda^{\nu+1}}\int_{0}^{\infty}\xi^{\nu}e^{-\xi
	}\left[  1-\frac{\xi^{\beta+1}}{(  \beta+1)  \lambda^{\beta+1}}+ \dots \right]  d\xi \qquad\text{as\ } \left\vert \lambda\right\vert \rightarrow\infty .
\end{equation}
Integration term by term gives 
\begin{equation}\label{e:supp:G:expand}
	G_{\beta,\nu}(  \lambda)  =\frac{\Gamma(\nu+1)}{\lambda^{\nu+1}}%
	-\frac{\Gamma(\nu+\beta+2)}{(  \beta+1)  \lambda^{\nu+\beta+2}%
	}+ \dots  \qquad \text{as\ }\left\vert \lambda\right\vert \rightarrow\infty . 
\end{equation}
Therewith, we obtain from~\eqref{supp:e:eval2} the approximated eigenvalue problem
\[
\vartheta \, \lambda=-G_{\beta,\nu}(
\lambda)  =-\frac{\Gamma(  \nu+1)  }{\lambda^{\nu+1}}+\frac
{\Gamma(  \nu+\beta+2)  }{(  \beta+1)  \lambda^{\nu+\beta+2}%
}+\dots \qquad\text{as\ }\left\vert \lambda\right\vert \rightarrow\infty .
\]
%Here, $G(0)$ depends on $\nu$, but it is always real positive and hence
%the specific value does not modify the stability of the solution. 
So, the next order expansion of the approximated eigenvalue equation~\eqref{e:supp:EV:zeta:asymp} is
\begin{equation}\label{e:supp:EV:zeta:second_order}
	\vartheta \lambda^{\nu+2}%
	\simeq -\Gamma(  \nu+1)  +\frac{\Gamma(  \nu+\beta+2)  }{(
		\beta+1)  \lambda^{\beta+1}} \qquad\text{for\ }\left\vert \lambda\right\vert \rightarrow\infty .
\end{equation}
From here, we find that the root with positive imaginary real part is approximated by
\begin{equation}\label{e:supp:zeta0}
	\lambda\simeq\lambda_0 = \biggl(\frac{\Gamma(  \nu+1)  }{\vartheta}\biggr)^{\frac{1}{\nu+2}}e^{\frac
		{i\pi}{\nu+2}}\qquad \text{as\ }\vartheta\rightarrow0
\end{equation}
Using the above formula in the second order expansion~\eqref{e:supp:EV:zeta:second_order} yields
\begin{equation}
	\vartheta \lambda^{\nu+2}%
	%=-\Gamma(  \nu+1)  +\frac{\Gamma(  \nu+\beta+2)  }{(\beta+1)  \lambda^{\beta+1}}
	\simeq-\Gamma(  \nu+1)  +\frac{\Gamma(\nu+\beta+2)}{(  \beta+1)  } \biggl(  \frac{\vartheta}{\Gamma(  \nu+1)  }\biggr)
	^{\frac{\beta+1}{\nu+2}}e^{-\frac{i\pi(  \beta+1)  }{\nu+2}}\quad\text{as\ } \vartheta\to 0. %
\end{equation}
From the above expansion, we deduce the next order correction $\tilde\lambda_2$ to $\lambda_0$ from~\eqref{e:supp:zeta0} by the expansion $\lambda\simeq \lambda_0 \bigl(1-\tilde\lambda_2\bigr)^{\frac{1}{\nu+2}}$ as $\vartheta\to 0$ with
\begin{equation}\label{e:supp:zeta0:2nd}
	\tilde\lambda_2 = \frac{\Gamma(  \nu+\beta+2)  }{(  \beta+1)
		\Gamma(  \nu+1)  } \biggl(  \frac{\vartheta}{\Gamma(  \nu+1)  }\biggr)
	^{\frac{\beta+1}{\nu+2}}e^{-\frac{i\pi(  \beta+1)  }{\nu+2}}.
\end{equation}
By using Taylor, we have $\lambda \simeq \lambda_0 \bigl( 1- \sfrac{\tilde\lambda_2}{\nu+2}\bigr)$.
% \[
% \lambda\simeq(  \frac{\Gamma(  \nu+1)  }{G(  0)
% (  \rfac)  ^{\frac{1}{\beta+1}}})  ^{\frac{1}{\nu+2}}e^{\frac
% {i\pi}{\nu+2}}\left[  1-\frac{\Gamma(  \nu+\beta+2)  }{(  \beta+1)
% (  \nu+2)  \Gamma(  \nu+1)  }(  \frac{G(
% 0)  (  \rfac)  ^{\frac{1}{\beta+1}}}{\Gamma(  \nu+1)
% })  ^{\frac{\beta+1}{\nu+2}}e^{-\frac{i\pi(  \beta+1)  }{\nu+2}}\right]
% \]
Next, we expand the exponential phase factor in~\eqref{e:supp:zeta0} as
\[
e^{\frac{i\pi}{\nu+2}}=\exp\biggl(\frac{i\pi}{2}+\frac{i\pi}{\nu+2}-\frac{i\pi}%
{2}\biggr)=i \exp\biggl(-\frac{i\pi \nu}{2(  \nu+2)  }\biggr)=i\biggl(  1-\frac{i\pi
	\nu}{2(  \nu+2)  }+ \dots \biggr)
\]
A combination of the expansion so far, leads to the asymptotic formula
\begin{equation}
	\lambda \simeq \biggl(\frac{\Gamma(  \nu+1)  }{\vartheta}\biggr)^{\frac{1}{\nu+2}} \bigl( i + \lambda_2 \bigr) \qquad\text{as\ } \rfac \to 0 
\end{equation}
with
\begin{equation}\label{e:supp:zeta2}
	\lambda_2 = \frac{\pi \nu}{2(  \nu+2)  }+\frac{\Gamma(  \nu+\beta+2)
	}{(  \beta+1)  (  \nu+2)  \Gamma(  \nu+1)  }
	\biggl(\frac{\vartheta}%
	{\Gamma(  \nu+1)  }\biggr)^{\frac{\beta+1}{\nu+2}} \; i \, e^{-\frac{i\pi(
			\beta+1)  }{\nu+2}} .
\end{equation}
% \[
% \lambda\simeq(  \frac{\Gamma(  \nu+1)  }{G(  0)
% (  \rfac)  ^{\frac{1}{\beta+1}}})  ^{\frac{1}{\nu+2}}i(
% 1-\frac{i\pi \nu}{2(  \nu+2)  }+...)  \left[  1-\frac
% {\Gamma(  \nu+\beta+2)  }{(  \beta+1)  (  \nu+2)
% \Gamma(  \nu+1)  }(  \frac{G(  0)  (
% \rfac)  ^{\frac{1}{\beta+1}}}{\Gamma(  \nu+1)  })
% ^{\frac{\beta+1}{\nu+2}}e^{-\frac{i\pi(  \beta+1)  }{\nu+2}}\right]
% \]%
% \[
% \lambda\simeq(  \frac{\Gamma(  \nu+1)  }{G(  0)
% (  \rfac)  ^{\frac{1}{\beta+1}}})  ^{\frac{1}{\nu+2}}i\left[
% 1-\frac{i\pi \nu}{2(  \nu+2)  }-\frac{\Gamma(  \nu+\beta+2)
% }{(  \beta+1)  (  \nu+2)  \Gamma(  \nu+1)  }(
% \frac{G(  0)  (  \rfac)  ^{\frac{1}{\beta+1}}}%
% {\Gamma(  \nu+1)  })  ^{\frac{\beta+1}{\nu+2}}e^{-\frac{i\pi(
% \beta+1)  }{\nu+2}}\right]
% \]
% \AS{I would prefer to write the above expansion as 
% \[
% \lambda\simeq(  \frac{\Gamma(  \nu+1)  }{G(  0)
% (  \rfac)  ^{\frac{1}{\beta+1}}})  ^{\frac{1}{\nu+2}}(i+\lambda_2) ,
% \]
Now, the real part of $\lambda_2$ determines the stability of the stationary state. By using $\nu \ll 1$, it is given by
% \[
%   \mathrm{Re}(\lambda_2) =  \frac{\pi \nu}{2(  \nu+2)  }-\frac{\Gamma(  \nu+\beta+2)
% }{(  \beta+1)  (  \nu+2)  \Gamma(  \nu+1)  }(
% \frac{G(  0)  (  \rfac)  ^{\frac{1}{\beta+1}}}%
% {\Gamma(  \nu+1)  })  ^{\frac{\beta+1}{\nu+2}} \sin(  \frac
% {\pi(  \beta+1)  }{\nu+2})
% \]
\begin{equation}\label{e:supp:null1:zeta2}
	\operatorname{Re}(\lambda_2)  \simeq \frac{\pi \nu}{4}-\frac{\Gamma(\beta+1)  }{2
	}\vartheta^{\frac{\beta+1}{\nu+2}}\sin\Bigl(
	\frac{\pi (  \beta+1)  }{2}\Bigr) \qquad\text{as }\ \nu\to 0 ,
\end{equation}
in accordance to~\eqref{e:null1:zeta2}.
% 
% The neutral curve separating stability from instability corresponds to the set
% of values of $\nu$ for which this root is purely imaginary. This is the curve:%
% \[
% -\frac{i\pi \nu}{2(  \nu+2)  }+\frac{i\Gamma(  \nu+\beta+2)
% }{(  \beta+1)  (  \nu+2)  \Gamma(  \nu+1)  }(
% \frac{G(  0)  (  \rfac)  ^{\frac{1}{\beta+1}}}%
% {\Gamma(  \nu+1)  })  ^{\frac{\beta+1}{\nu+2}}\sin(  \frac
% {\pi(  \beta+1)  }{\nu+2})  =0
% \]
% Equivalently this is the curve:%
% \[
% \frac{\pi \nu}{2(  \nu+2)  }=\frac{\Gamma(  \nu+\beta+2)
% }{(  \beta+1)  (  \nu+2)  \Gamma(  \nu+1)  }(
% \frac{G(  0)  (  \rfac)  ^{\frac{1}{\beta+1}}}%
% {\Gamma(  \nu+1)  })  ^{\frac{\beta+1}{\nu+2}}\sin(  \frac
% {\pi(  \beta+1)  }{\nu+2})
% \]
% where $n\rightarrow0$ and $\beta$ is in principle arbitrary. Approximating this
% equation to the leading order in $\nu$ we obtain:%
% \[
% \frac{\pi \nu}{4}\simeq\frac{\Gamma(  \beta+1)  }{2
% }G(  0)^{\frac{\beta+1}{\nu+2}}  (  \rfac)  ^{\frac{1}{\nu+2}}\sin(
% \frac{\pi(  \beta+1)  }{2})
% \]
% Notice that there is a dependence of $\nu$ in $\rfac.$ Therefore there is not
% just a relation between $\nu$ and $\beta.$ For each $\nu>0$ we obtain stability if
% $\rfac\rightarrow0.$ On the other hand we obtain a curve of stability with
% $\nu$ depending on $\rfac$ if $\sin(  \frac{\pi(  \beta+1)  }%
% {2})  >0.$
 % arXiv

\bibliographystyle{siamplain}
\bibliography{references}
\end{document}